\documentclass[11pt]{article}
\small
\normalsize
\usepackage{amsmath}
\usepackage{amssymb}
\usepackage{amsfonts}
\usepackage{amsbsy}
\usepackage{amscd}
\usepackage{theorem}
\usepackage{epsf}
\usepackage{psfrag}
\usepackage{multicol}
\allowdisplaybreaks
\font\goth=eufm10 scaled 1200

\def\noi{{\noindent}}
\def\cq{ \hfill $\blacksquare$ }
\def\llbracket{[\hspace{-.10em} [ }
\def\rrbracket{ ] \hspace{-.10em}]}

\def\s{{\cal S}}

\def\p{{\cal P}}
\def\m{{\cal M}}
\def\l{{\cal L}}
\def\h{{\cal H}}
\def\ii{{\cal I}}
\def\jj{{\cal J}}
\def\e{{\cal E}}
\def\t{{\cal T}}

\def\z{{\cal Z}}
\def\n{{\cal N}}

\def\E{{\bf E}}
\def\P{{\bf P}}

\def\ov{\overline}

\def\wt{\widetilde}

\def\la{\longrightarrow}
\def\da{\downarrow}

\def\dimh{{\rm dim}_h}
\def\dimp{{\rm dim}_p}
\def\diml{\underline{ {\rm dim} } }
\def\dimu{\overline{ {\rm dim} } }

\def\noi{\noindent}

\def\proof{\noindent{\bf Proof.} }

\def\supp{{\rm supp\,}}
\def\build#1_#2^#3{\mathrel{\mathop{\kern 0pt#1}\limits_{#2}^{#3}}}

\newtheorem{theorem}{Theorem}[section]
\newtheorem{lemma}[theorem]{Lemma}
\newtheorem{proposition}[theorem]{Proposition}
\newtheorem{corollary}[theorem]{Corollary}
\newtheorem{definition}{Definition}[section]

\newcommand{\R}{\mathbb{R}}
\newcommand{\Z}{\mathbb{Z}}

\newcommand{\Q}{\mathbb{Q}}

\newcommand{\N}{\mathbb{N}}
\newcommand{\M}{\mathbb{M}}
\newcommand{\T}{\mathbb{T}}

\newcommand{\un}{\boldsymbol{1}}

\begin{document}

\title{ \bf PROBABILISTIC AND FRACTAL ASPECTS \\
OF LEVY TREES} 
\author{ by \\
Thomas {\sc Duquesne,}\\
{\small Universit\'e Paris 11, Math\'ematiques, 91405 Orsay Cedex, France } \\
and  \\
Jean-Fran\c cois {\sc Le Gall} \\
{\small D.M.A., Ecole normale sup\'erieure, 45 rue d'Ulm, 75005 Paris, France} }
\vspace{4mm}
\date{\today} 

\maketitle

\begin{abstract} We investigate the random continuous trees called L\'evy trees, which 
are obtained as scaling limits of discrete Galton-Watson trees.
We give a mathematically precise definition of these random trees as random variables taking
values in the set of equivalence classes of compact rooted $\R$-trees, which is equipped
with the Gromov-Hausdorff distance. To construct L\'evy trees, we make
use of the coding by the height process which was studied in detail in previous work. We then investigate
various probabilistic properties of L\'evy trees. In particular we establish a branching property
analogous to the well-known property for Galton-Watson trees: Conditionally given the tree
below level $a$, the subtrees originating from that level are distributed as the atoms of
a Poisson point measure whose intensity involves a local time measure supported on the 
vertices at distance $a$ from the root. We study regularity properties of local times
in the space variable, and prove that the support of local time is the full level set, except
for certain exceptional values of $a$ corresponding to local extinctions. We also compute several
fractal dimensions of L\'evy trees, including Hausdorff and packing dimensions, in terms
of lower and upper indices for the branching mechanism function $\psi$ which characterizes the distribution
of the tree. We finally discuss some applications to super-Brownian motion with a general branching
mechanism.
\end{abstract}


\section{Introduction.}

This work is devoted to the study of various properties of the so-called L\'evy trees, which
are continuous analogues of the discrete Galton-Watson trees. Our main contributions to the
probabilistic analysis of L\'evy trees include the construction of local time measures
supported on level sets of the tree, the use of these local times to formulate and
establish a branching property analogous to a well-known result in the discrete setting,
and the proof of a ``subtree'' decomposition along the ancestral line of a typical
vertex in the tree. Additionally, we study the
fractal properties of L\'evy trees and compute their Hausdorff and packing dimensions as well as that of
particular subsets such as level sets, under broad assumptions on the branching mechanism
characterizing the tree.

One major originality of the present article compared to our previous work
\cite{DuLG},\cite{LGLJ1},\cite{LGLJ2} is to view L\'evy trees as random variables
taking values in the space of compact rooted $\R$-trees. The precise definition of
an $\R$-tree is recalled in Section 2 below. Informally 
an $\R$-tree is a metric space $(\t,d)$ such that for any two points $\sigma$ and $\sigma'$
in $\t$ there is a unique arc with endpoints $\sigma$ and $\sigma'$ and furthermore this
arc is isometric to a compact interval of the real line. A rooted $\R$-tree 
is an $\R$-tree with a distinguished vertex called the root. We write $h(\t)$
for the height of $\t$, that is the maximal distance from the root to a vertex in $\t$. 
Say that two
rooted $\R$-trees are equivalent if there is a root-preserving isometry
that maps one onto the other. It was noted in
\cite{EPW} that the set of equivalence classes of compact rooted $\R$-trees, equipped 
with the Gromov-Hausdorff distance \cite{Gro} is a Polish space. 

The study of $\R$-trees has been motivated by algebraic and geometric purposes.
See in particular \cite{Pau} and the survey \cite{DMT}. One of our goals is
to initiate a probabilistic theory of $\R$-trees, by starting with the fundamental
case of L\'evy trees. See \cite{EPW} for another probabilistic application
of $\R$-trees. We also mention the 
recent article \cite{AMP}, which discusses a different class of continuum random trees
obtained as weak limits of birthday trees (instead of the Galton-Watson trees
considered here), using ideas 
related to the present work.

To motivate our definition of L\'evy trees, let us describe a simple approximation
result, which is a special case of Theorem \ref{discrete-continuous} below. Let $\mu$
be a probability measure on $\Z_+$, with $\mu(1)<1$. Assume that $\mu$ has mean one and is
in the domain of attraction of a stable distribution with index $\gamma\in(1,2]$. 
When $\gamma=2$, this holds as soon as $\mu$ has finite variance, and when $\gamma\in(1,2)$,
it is enough to assume that $\mu(k)\sim c\,k^{-1-\gamma}$ as $k\to\infty$. Denote by 
$\theta$ a Galton-Watson tree with offspring distribution $\mu$, which describes the
genealogy of a (discrete-time) Galton-Watson branching process with offspring distribution $\mu$
started initially with one ancestor. We can view $\theta$ as a (random) finite graph
and equip it with the natural graph distance. If $r>0$, the scaled tree 
$r\theta$ is obviously defined by requiring the distance between two neighboring
vertices to be $r$ instead of $1$. Also let $h(\theta)$ stand for
the  maximal generation in $\theta$. Then there is a $\sigma$-finite measure $\Theta(d\t)$
on the space of (equivalence classes of) rooted compact $\R$-trees such that for every
$a>0$, the conditional law of the scaled tree $n^{-1}\theta$ knowing that $h(\theta)\geq an$
converges as $n\to\infty$ to the probability measure $\Theta(d\t\mid h(\t)\geq a)$,
in the sense of weak convergence for the Gromov-Hausdorff distance on pointed metric spaces.

In a sense, the preceding result is not really new: See \cite{Al2},\cite{Duq} and especially
Chapter 2 of \cite{DuLG} for related limit theorems with a different formalism. Still we believe 
that the formalism of $\R$-trees is useful both to formulate such results and
to analyse the limiting objects as we do in the present work. 

Let us turn to a more precise description of the class of
random trees that will be considered here. A L\'evy tree can be interpreted
as the genealogical tree of a continuous-state branching process, whose law 
is characterized by a real function $\psi$ defined on $[0,\infty)$,
which is called the branching
mechanism. Here we restrict our attention to the critical or subcritical case where 
$\psi$ is nonnegative and of the form
$$ \psi (\lambda )= \alpha \lambda + \beta \lambda^2 + 
\int_{(0, \infty)} \pi (dr) (e^{-\lambda r} -1+\lambda r) \; , \quad \lambda\geq 0,  $$
where $\alpha, \beta \geq 0$ and $\pi $ is a $\sigma $-finite measure on $(0, \infty)$ such that 
$\int_{(0, \infty)} \pi (dr)(r\wedge r^2) < \infty $. 
We assume throughout this work the condition
$$\int_1^\infty {du\over\psi(u)}<\infty$$
which is equivalent to the a.s. extinction of the 
continuous-state branching process, and thus necessary for the compactness
of the associated genealogical tree. Of particular importance are the
quadratic branching case $\psi(\lambda)=c\,\lambda^2$ and the stable case
$\psi(\lambda)=c\,\lambda^\gamma$, $1<\gamma<2$, which both arise in the discrete
approximation described above.

The precise definition of the $\psi$-L\'evy tree then
depends on the height process introduced by Le Gall and Le Jan \cite{LGLJ1}
(see also Chapter 1 of \cite{DuLG}) in view of coding the
genealogy of general continuous-state branching processes.
The height process is obtained as a functional of the spectrally positive L\'evy process $X$
with Laplace exponent $\psi$. An important role is played by the excursion measure $N$ of $X$
above its minimum process.
In the quadratic branching case $\psi(u)=c\,u^2$, $X$  is a (scaled) Brownian motion, the height process
$H$ is
a reflected Brownian motion and the ``law'' of $H$ under $N$ is just the
It\^o measure of positive excursions of linear Brownian motion: This is related to the fact that
the contour  process of Aldous' {\it Continuum Random Tree} is given by a normalized
Brownian excursion (see \cite{Al1} and \cite{Al2}), or to the Brownian snake construction
of superprocesses with quadratic branching mechanism (see e.g. \cite{LG1}).
In our more general setting, the height  process can be defined informally as
follows. For every $t\geq 0$, $H_t$ measures
the size of the set $\{s\leq t:X_s=\inf_{[s,t]}X_r\}$. A precise definition of $H_t$
is recalled in Section 3 below. Under our assumptions, the process $H$ has a continuous
modification. 

The claim is now that the sample path of $H$ under $N$ codes a random continuous
tree called the $\psi$-L\'evy tree. The precise meaning of the coding is explained
in Section 2 in a deterministic setting, but let us immediately outline the construction
of the tree. 
We write $\zeta$ for the duration of the excursion under $N$ and define 
a random function $d_H$ on $[0,\zeta]^2$ by setting
$$ d_H (s,t) = H_s + H_t -2m_H (s,t) \; ,$$
where we have set $ m_H (s,t)= \inf_{s\wedge t \leq r\leq s\vee t} H_r $. 
We introduce an associated equivalence relation by setting $s\thicksim_H t$ if and only if $d_H(s,t)=0$. 
In particular, $0\thicksim_H \zeta$. 
The function $d_H$ obviously extends to the quotient set $\t_H:=[0, \zeta] / \thicksim_H $ 
and defines a distance on this set. It is not hard to verify that 
$(\t_H,d_H)$ is a compact $\R$-tree, and its root is by definition the equivalence
class of $0$. 
Informally, each real number $s\in[0,\zeta]$ corresponds to a vertex at level $H_s$ in the
tree, and $d_H(s,t)$ is the distance between vertices corresponding to $s$ and $t$
(in particular $s$ and $t$ correspond to the same vertex if and only if $d_H(s,t)=0$). The 
quantity $m_H(s,t)$
can be interpreted as the generation of the most recent common ancestor to $s$ and $t$.

The law of the L\'evy tree is by
definition the distribution $\Theta(d\t)$ of the compact rooted $\R$-tree
$(\t_H,d_H)$ under the measure $N$. Notice that $N$ is an infinite measure, and so is $\Theta$. 
However, for every $a>0$, $v(a):=\Theta(h(\t)>a)<\infty$, and more precisely $v(a)$ is determined
by the equation
$$\int_{v(a)}^\infty {du\over \psi(u)}=a.$$ 
Section 4 contains the proof of several important properties of L\'evy trees. 
In particular, for every $a> 0$, we construct the local time $\ell^a$
at level $a$, which is a finite measure supported on the level set
$$\t(a):=\{\sigma\in\t:d(\rho(\t),\sigma)=a\}$$ 
where $\rho(\t)$ denotes the root of $\t$. We then prove the
fundamental ``branching property'': If $(\t^{(i),\circ},i\in \ii)$ denote the connected components
of the open set $\{\sigma\in\t:d(\rho(\t),\sigma)>a\}$, the closure 
$\t^{(i)}$ of each $\t^{(i),\circ}$ is a compact rooted $\R$-tree with root
$\sigma_i\in\t(a)$ and, conditionally on $\ell^a$, the point measure
$$\sum_{i\in \ii} \delta_{(\sigma_i,\t^{(i)})}$$
is Poisson with intensity $\ell^a(d\sigma)\,\Theta(d\t)$ (see Theorem \ref{existLT}
for a slightly more precise result stating that this point measure is also
independent of the part of the tree ``below level $a$''). Up to some point, the
branching property follows from a result of
\cite{DuLG} (Proposition 1.3.1 or Proposition 4.2.3) showing that excursions of the height process above
level $a$ are distributed as the atoms of a Poisson point measure whose intensity 
is (a random multiple of) the law of $H$ under $N$.
In this form, the branching property has been recently used by Miermont \cite{Mie} to investigate
self-similar fragmentations of the stable tree.

Using the branching property, we investigate the regularity properties of 
local times. We show that the mapping $a\la \ell^a$ has a c\` adl\` ag modification
and that, except for a countable set of values of $a$ (corresponding to
local extinctions of the tree) the support of $\ell^a$ is 
the full level set $\t(a)$. This is used in Section 6 to extend to
superprocesses with a general branching mechanism a continuity property of
the support process that had been derived by Perkins \cite{Per} in the 
quadratic case.

In the final part of Section 4 we prove a Palm-like decomposition of the tree
along the ancestor line of a typical vertex at level $a$ (Theorem \ref{Palmdec}). This decomposition plays
an important role in Section 5. We use it in Section 4 to analyse the multiplicity
of vertices of the tree. By definition, the multiplicity $n(\sigma)$ of $\sigma\in\t$ is the
number of connected components  of $\t\backslash\{\sigma\}$. We prove that 
$\Theta$ a.e. $n(\sigma)$ takes values in the set $\{1,2,3,\infty\}$. We also
characterize the branching mechanism functions $\psi$ for which there exist
binary ($n(\sigma)=3$) or infinite ($n(\sigma)=\infty$) branching points.
We then observe that infinite branching points are related to discontinuities of
local times: Precisely, for any level $b$ such that the mapping $a\longrightarrow \ell^a$
is discontinuous at $b$, there is a (unique) infinite branching point $\sigma_b$
such that $\ell^b=\ell^{b-}+\lambda_b\,\delta_{\sigma_b}$ for some $\lambda_b>0$.
As a last application of our Palm decomposition, we prove an invariance property
of the measure $\Theta$ under uniform re-rooting (Proposition \ref{re-rooting}).
   
Section 5 is mostly devoted to the computation of the Hausdorff 
and packing dimensions of various subsets of
$\t$.  For any subset $A$ of $\t$, we denote by $\dimh (A)$ the Hausdorff dimension of $A$ and by $\dimp
(A)$ its packing dimension.  Following \cite{Fal}, Section 3.1, we also consider the
lower and upper box counting dimensions of $A$:
$$ \diml (A)= \liminf_{\delta \rightarrow 0} \frac
{\log \left( \n (A, \delta ) \right)}{\log (1/\delta)} \ ,\quad 
\dimu (A)= \limsup_{\delta \rightarrow 0} \frac{\log \left( \n (A, \delta ) \right)}{\log (1/\delta)}\ ,$$
where $\n (A, \delta )$ is the minimal number of open balls with
radius
$\delta $ that are necessary to cover $A$.
In order to state our main results, we need to 
introduce the lower and upper indices of $\psi $ at infinity:
$$\gamma = \sup \{ a\geq 0 \; :\; \lim_{\lambda \rightarrow \infty } 
\lambda^{-a}\psi (\lambda ) =+\infty  \} \quad , 
\quad \eta =\inf \{ a\geq 0 \; :\; \lim_{\lambda \rightarrow \infty } 
\lambda^{-a}\psi (\lambda ) =0  \}. $$
Note that $1\leq \gamma \leq \eta $ and that $\eta=\gamma$ if $\psi$ is 
regularly varying at infinity. 
Let $E$ be a nonempty compact subset of the interval $(0,\infty)$ and assume that $E$ is regular in the
sense that its Hausdorff and  upper box counting dimensions coincide: $\dimh(E)=\dimu(E) =d(E) \in [0, 1] $
(here
$\dimh$ and $\dimu$ obviously refer to the usual metric on the real line). Set $a=\sup E$ and 
$$ \t(E)= \bigcup_{l\in E} \t(l).$$
Theorem \ref{hauspacktE} asserts that under the assumption $\gamma > 1$, we have $\Theta$-a.e. 
on $\{ h(\t)>a\}$,  
$$ \diml (\t(E) ) =\dimh (\t(E) ) = d(E) + \frac{1}{\eta -1} \quad {\rm and }  \quad 
\dimu (\t(E) ) =\dimp (\t(E) ) =  d(E) + \frac{1}{\gamma -1}.$$
In particular, 
we have $\Theta$-a.e. 
$$\dimh (\t) = \frac{\eta}{\eta -1}\quad,\quad\dimp (\t ) = \frac{\gamma }{\gamma -1}$$
and, $\Theta$-a.e. on $\{ h(\t)>a \}$,
$$\dimh (\t(a)) = \frac{1}{\eta -1}\quad,\quad\dimp (\t(a) ) = \frac{1}{\gamma -1}.$$
Note that in the stable branching case $\psi(u)=u^\gamma$, the Hausdorff dimension of
$\t$ has been computed by Haas and Miermont \cite{HaM} independently of the present work.

The proofs rely on the classical results linking upper and lower densities of a measure with 
the Hausdorff and packing dimensions of its support.
Another useful ingredient is the following estimate for covering numbers  of $\t $
(Proposition \ref{covt}). We have $\Theta$-a.e. for all 
sufficiently small $\delta $, 
$$ \frac{v(2\delta )}{4\delta }\,\zeta \leq \n (\t , \delta ) \leq \frac{12v(\delta /6)}{\delta }\,\zeta .$$

In Section 6, we give an application of Theorem \ref{hauspacktE} 
to the range 
of a superprocess 
$Z =(Z_l ,  l\geq 0)$ with branching mechanism $\psi $, whose 
spatial motion is standard Brownian motion in $\R^k$. To this end, we 
introduce the notion of a spatial tree, which allows us to combine
the genealogical structure of $\t$ with independent spatial Brownian motions.
Of course, this is more or less equivalent to the L\'evy snake approach 
of \cite{LGLJ2} and Chapter 4 of \cite{DuLG}. Still the formalism of
$\R$-trees makes this construction more tractable and more efficient for
applications. Roughly speaking, spatial trees allow us to express the superprocess $Z$ in terms of 
the occupation measure of a Gaussian process indexed by $\t$. It is therefore possible 
to use soft arguments to lift
fractal  properties of the index set $\t$ to the range of $Z$.  We prove the following result (Theorem
\ref{Haus-super}).  Let $E\subset (0,\infty)$ and $a=\sup E$ be as above. Denote by $R_E $ the range of
$Z$ over the time set
$E$, defined by
$$R_E= 
\overline{ \bigcup_{l\in E} {\rm supp }\, Z_l  } $$
where ${\rm supp }\, Z_l  $ stands for the topological support of $Z_l$. If $\gamma >1$, then a.s. on 
$\{ \langle  Z_a, 1 \rangle \neq 0 \}$,
\begin{equation}
\label{supdim}
\dimh (R_E )= \left( 2d(E) +  \frac{2}{\eta -1} \right)\wedge k.
\end{equation}
In the quadratic branching case, this result was obtained earlier by
Tribe \cite{Tribe} (see also Serlet \cite{Ser95}). For more general superprocesses,
closely related results can be found in Theorem 2.1 of Delmas \cite{Del}, whose proof depends
on a subordination method which requires certain restrictions on the 
branching mechanism function $\psi$. See Dawson \cite{Daw} for more references and
results in the stable branching case.

This paper is intended to be as self-contained as possible. However, it is clear 
that many of our results depend on properties of the height process $H$ that were derived in the
monograph \cite{DuLG}. For the reader's convenience, we have recalled most of
the needed results in Section 3 below.

The paper is organized as follows. Section 2
explains the coding of trees by continuous functions 
in a deterministic setting, and also includes a brief
discussion of the convergence of trees in the Gromov-Hausdorff
metric. 
Section 3 recalls the basic facts about the height process and 
establishes an important preliminary result that is needed for the ancestral line decomposition of 
subsection 4.3. Section 4 is the core of this paper. It first contains
the precise definition of the L\'evy
tree as the tree
coded by the excursion of $H$ under $N$, in the framework of Section 2. This definition is
justified by limit theorems relating discrete and continuous trees. Section
4 then presents the basic probabilitic properties of L\'evy trees, in particular
the branching property, the existence and regularity of local times and the 
decomposition along an ancestral line. 
Fractal properties of L\'evy trees are studied in
Section 5. Finally, Section 6  discusses applications to superprocesses. 

\section{Deterministic trees}

\subsection{The $\R$-tree coded by a continuous function}

We start with a basic definition (see e.g. \cite{DMT}).

\begin{definition}
A metric space $(\t,d)$ is an $\R$-tree if the following two
properties hold for every $\sigma_1,\sigma_2\in \t$.

\begin{description}
\item{(i)} There is a unique
isometric map
$f_{\sigma_1,\sigma_2}$ from $[0,d(\sigma_1,\sigma_2)]$ into $\t$ such
that $f_{\sigma_1,\sigma_2}(0)=\sigma_1$ and $f_{\sigma_1,\sigma_2}(
d(\sigma_1,\sigma_2))=\sigma_2$.
\item{(ii)} If $q$ is a continuous injective map from $[0,1]$ into
$\t$, such that $q(0)=\sigma_1$ and $q(1)=\sigma_2$, we have
$$q([0,1])=f_{\sigma_1,\sigma_2}([0,d(\sigma_1,\sigma_2)]).$$
\end{description}

\noindent A rooted $\R$-tree is an $\R$-tree $(\t,d)$
with a distinguished vertex $\rho=\rho(\t)$ called the root.
\end{definition}

In what follows, $\R$-trees will always be rooted, even if this
is not mentioned explicitly. 

Let us consider a rooted $\R$-tree $(\t,d)$.
The range of the mapping $f_{\sigma_1,\sigma_2}$ in (i) is denoted by
$\llbracket \sigma_1,\sigma_2\rrbracket$ (this is the line segment between $\sigma_1$
and $\sigma_2$ in the tree). 
In particular, for every $\sigma\in \t$, $\llbracket \rho,\sigma\rrbracket$ is the path 
going from the root to $\sigma$, which we will interpret as the ancestral
line of vertex $\sigma$. More precisely we define a partial order on the
tree by setting $\sigma\preccurlyeq \sigma'$
($\sigma$ is an ancestor of $\sigma'$) if and only if $\sigma\in\llbracket \rho,\sigma'
\rrbracket$.

If $\sigma,\sigma'\in\t$, there is a unique $\eta\in\t$ such that
$\llbracket \rho,\sigma
\rrbracket\cap \llbracket \rho,\sigma'
\rrbracket=\llbracket \rho,\eta
\rrbracket$. We write $\eta=\sigma\wedge \sigma'$ and call $\eta$ the most recent
common ancestor to $\sigma$ and $\sigma'$.

By definition, the multiplicity of a vertex $\sigma\in\t$ is the number of
connected components of $\t\backslash \{\sigma\}$. Vertices
of $\t\backslash\{\rho\}$ which have multiplicity
$1$ are called leaves.

\medskip

Our main goal in this section is to describe a method for constructing $\R$-trees, which
is particularly well-suited to our forthcoming applications to random trees.
We consider a	 (deterministic) continuous function
$g:[0,\infty)\longrightarrow[0,\infty)$ with compact support
and such that $g(0)=0$. To avoid
trivialities, we will also assume that $g$ is not identically zero.
For every $s,t\geq 0$, we set
$$m_g(s,t)=\inf_{r\in[s\wedge t,s\vee t]}g(r),$$
and
$$d_g(s,t)=g(s)+g(t)-2m_g(s,t).$$
Clearly $d_g(s,t)=d_g(t,s)$ and it is also easy to verify the triangle
inequality
$$d_g(s,u)\leq d_g(s,t)+d_g(t,u)$$
for every $s,t,u\geq 0$. We then introduce the equivalence relation
$s\sim t$ iff $d_g(s,t)=0$ (or equivalently iff $g(s)=g(t)=m_g(s,t)$). Let
$\t_g$ be the quotient space
$$\t_g=[0,\infty)/ \sim.$$
Obviously the function $d_g$ induces a distance on $\t_g$, and we keep the
notation $d_g$ for this distance. We denote by
$p_g:[0,\infty)\longrightarrow
\t_g$ the canonical projection. Clearly $p_g$ is continuous (when
$[0,\infty)$ is equipped with the Euclidean metric and $\t_g$ with the
metric $d_g$).

\begin{theorem}
\label{tree-deterministic}
The metric space $(\t_g,d_g)$ is an $\R$-tree.
\end{theorem}

We will view $(\t_g,d_g)$ as a rooted $\R$-tree with root $\rho=p_g(0)$.
If $\zeta>0$ is the supremum of the support of $g$, we have
$p_g(t)=\rho$ for every 
$t\geq \zeta$. In particular, $\t_g=p_g([0,\zeta])$ is compact.
We will call $\t_g$ {\it the $\R$-tree coded by $g$}.

Before proceeding to the proof of the theorem, we state and prove the
following root change lemma. 

\begin{lemma}
\label{root-change}
Let $s_0\in[0,\zeta)$. For any real $r\geq 0$, denote by
$\ov r$ the unique element of $[0,\zeta)$ such that
$r-\ov r$ is an integer multiple of $\zeta$. Set
$$g'(s)=g(s_0)+g(\ov{s_0+s})-2m_g(s_0,\ov{s_0+s}),$$
for every $s\in [0,\zeta]$, and $g'(s)=0$ for $s>\zeta$. Then,
the function $g'$ is continuous with compact support and satisfies
$g'(0)=0$, so that we can define
the metric space $(\t_{g'},d_{g'})$. Furthermore, for
every
$s,t\in[0,\zeta]$, we have
\begin{equation}
\label{iso-root}
d_{g'}(s,t)=d_g(\ov{s_0+s},\ov{s_0+t})
\end{equation}
and there exists a unique isometry $R$ from $\t_{g'}$ onto $\t_{g}$
such that, for every $s\in [0,\zeta]$,
\begin{equation}
\label{def-changeroot}
R(p_{g'}(s))=p_g(\ov{s_0+s}).
\end{equation}
\end{lemma}

\noindent{\bf Proof.} It is immediately checked that $g'$ satisfies
the same assumptions as $g$, so that we
can make sense of the tree $\t_{g'}$. Then the key step is to verify
the relation (\ref{iso-root}). Consider first the case where
$s,t\in[0,\zeta-s_0)$. Then two possibilities may occur.

If $m_g(s_0+s,s_0+t)\geq m_g(s_0,s_0+s)$, then
$m_g(s_0,s_0+r)=m_g(s_0,s_0+s)=m_g(s_0,s_0+t)$ for every $r\in[s,t]$, and
so
$$m_{g'}(s,t)=g(s_0)+m_g(s_0+s,s_0+t)
-2m_g(s_0,s_0+s).$$
It follows that
\begin{eqnarray*}
d_{g'}(s,t)&=&g'(s)+g'(t)-2m_{g'}(s,t)\\
&=&g(s_0+s)-2m_g(s_0,s_0+s)+g(s_0+t)-2m_g(s_0,s_0+t)\\
&&\qquad-2(m_g(s_0+s,s_0+t)-2m_g(s_0,s_0+s))\\
&=&g(s_0+s)+g(s_0+t)
-2m_g(s_0+s,s_0+t)\\
&=&d_g(s_0+s,s_0+t).
\end{eqnarray*}

If $m_g(s_0+s,s_0+t)< m_g(s_0,s_0+s)$, then the minimum in the definition
of $m_{g'}(s,t)$ is attained at $r_1$ defined as the first $r\in[s,t]$
such that $g(s_0+r)=m_g(s_0,s_0+s)$ (because for $r\in[r_1,t]$ we will
have $g(s_0+r)-2m_g(s_0,s_0+r)\geq -m_g(s_0,s_0+r)\geq
-m_g(s_0,s_0+r_1)$). Therefore,
$$m_{g'}(s,t)=g(s_0)-m_g(s_0,s_0+s),$$
and
\begin{eqnarray*}
d_{g'}(s,t)
&=&g(s_0+s)-2m_g(s_0,s_0+s)+g(s_0+t)-2m_g(s_0,s_0+t)
+2m_g(s_0,s_0+s)\\
&=&d_g(s_0+s,s_0+t).
\end{eqnarray*}

The other cases are treated in a similar way and are left to the reader.

By (\ref{iso-root}), if $s,t\in[0,\zeta]$ are such that $d_{g'}(s,t)=0$,
we have
$d_g(\ov{s_0+s},\ov{s_0+t})=0$ so that $p_g(\ov{s_0+s})=p_g(\ov{s_0+t})$.
Noting that $\t_{g'}=p_{g'}([0,\zeta])$ (the supremum of the support of
$g'$
is less than or equal to $\zeta$), we can define
$R$ in a unique way by the relation (\ref{def-changeroot}). From
(\ref{iso-root}),
$R$ is an isometry, and it is also immediate that
$R$ takes $\t_{g'}$ onto $\t_g$. \cq

\medskip
\noindent{\bf Proof of Theorem \ref{tree-deterministic}.}
Let us start with some 
preliminaries. For $\sigma,\sigma'\in\t_g$, we set
$\sigma\preccurlyeq\sigma'$ if and only if $d_g(\sigma,\sigma')
=d_g(\rho,\sigma')-d_g(\rho,\sigma)$.
If $\sigma=p_g(s)$ and $\sigma'=p_g(t)$, it follows from our definitions
that $\sigma\preccurlyeq\sigma'$ iff $m_g(s,t)=g(s)$.
It is immediate to verify that this defines a partial order on $\t_g$.

For any $\sigma_0,\sigma\in \t_g$, we set
$$\llbracket \sigma_0,\sigma\rrbracket=\{\sigma'\in
\t_g:d_g(\sigma_0,\sigma)=d_g(\sigma_0,\sigma')+d_g(\sigma',\sigma)\}.$$
If $\sigma=p_g(s)$ and $\sigma'=p_g(t)$,
then it is easy to verify that $\llbracket\rho,\sigma\rrbracket
\cap\llbracket\rho,\sigma'\rrbracket=\llbracket\rho,
\gamma\rrbracket$, where $\gamma=p_g(r)$, if $r$ is any time
which achieves the minimum of $g$ between $s$ and $t$. We then put $\gamma=\sigma\wedge
\sigma'$.

We set
$\t_g[\sigma]:=\{\sigma'\in\t_g:\sigma\preccurlyeq\sigma'\}$. If
$\t_g[\sigma]\not =\{\sigma\}$ and $\sigma\not =\rho$, then
$\t_g\backslash\t_g[\sigma]$ and $\t_g[\sigma]\backslash\{\sigma\}$
are two nonempty disjoint open sets. To see that $\t_g\backslash
\t_g[\sigma]$ is open, let $s$ be such that $p_g(s)=\sigma$ and note that
$\t_g[\sigma]$ is the image under  $p_g$ of the compact set
$\{u\in[0,\zeta]:m_g(s,u)=g(s)\}$. The set
$\t_g[\sigma]\backslash\{\sigma\}$ is open because if $\sigma'\in
\t_g[\sigma]$ and $\sigma'\ne \sigma$, it easily follows from our definitions that the open ball
centered at $\sigma'$ with radius $d_g(\sigma,\sigma')$ is contained
in $\t_g[\sigma]\backslash\{\sigma\}$.

We now prove property (i)
of the definition of an $\R$-tree. By using Lemma \ref{root-change}
with $s_0$ such that $p_g(s_0)=\sigma_1$, we may assume that
$\sigma_1=\rho=p_g(0)$. If $\sigma\in\t_g$ is fixed, we have to prove that
there exists a unique isometry $f$ from $[0,d_g(\rho,\sigma)]$ into
$\t_g$ such that $f(0)=\rho$ and
$f(d_g(\rho,\sigma))=\sigma$. Let $s\in p_g^{-1}(\{\sigma\})$, so
that $g(s)=d_g(\rho,\sigma)$. Then, for every
$a\in[0,d_g(\rho,\sigma)]$, we set
$$w(a)=\inf\{r\in[0,s]:m_g(r,s)=a\}.$$
Note that $g(w(a))=a$. We put
$f(a)=p_g(w(a))$. We have $f(0)=\rho$ and
$f(d_g(\rho,\sigma))=\sigma$, the latter because
$m_g(w(g(s)),s)=g(s)$ implies $p_g(w(g(s)))=p_g(s)=\sigma$. It is also
easy to verify that
$f$ is an isometry: If $a,b\in[0,d_g(\rho,\sigma)]$
with $a\leq b$, it is immediate that $m_g(w(a),w(b))=a$, and so
$$d_g(f(a),f(b))=g(w(a))+g(w(b))-2a=b-a.$$

To get uniqueness, suppose that $\tilde f$ is an isometry satisfying the
property in (i). Then, if $a\in[0,d_g(\rho,\sigma)]$,
$$d_g(\tilde f(a),\sigma)=d_g(\rho,\sigma)-a=d_g(\rho,\sigma)
-d_g(\rho,\tilde f(a)).$$
Therefore, $\tilde f(a)\preccurlyeq \sigma$. Recall that
$\sigma=p_g(s)$, and choose $t$ such that $p_g(t)=\tilde f(a)$.
Note that $g(t)=d_g(\rho,p_g(t))=a$. Since
$\tilde f(a)\preccurlyeq \sigma$ we have $g(t)=m_g(t,s)$.
On the other hand, we also know that $a=g(w(a))=m_g(w(a),s)$.
It follows that we have $a=g(t)=g(w(a))=m_g(w(a),t)$
and thus $d_g(t,w(a))=0$, so that
$\tilde f(a)=p_g(t)=p_g(w(a))=f(a)$. This completes the proof of (i).

As a by-product of the preceding argument, we see that
$f([0,d_g(\rho,\sigma)]) =\llbracket\rho,\sigma\rrbracket$:
Indeed, we have seen that for every
$a\in[0,d_g(\rho,\sigma)]$, we have $f(a)\preccurlyeq \sigma$
and, on the other hand, if $\eta\preccurlyeq \sigma$, the end of the
proof of (i) just shows that $\eta=f(d_g(\rho,\eta))$.

We turn to the proof of (ii). We let $q$ be a continuous injective
mapping from $[0,1]$ into $\t_g$, and we aim at proving
that $q([0,1])=f_{q(0),q(1)}([0,d_g(q(0),q(1))])$. From Lemma
\ref{root-change} again, we may assume that $q(0)=\rho$,
and we set $\sigma=q(1)$. Then we have just noticed that
$f_{0,\sigma}([0,d_g(\rho,\sigma)])=\llbracket \rho,\sigma
\rrbracket$.

We first argue by contradiction to prove that
$\llbracket \rho,\sigma
\rrbracket\subset q([0,1])$. Suppose that
$\eta\in \llbracket \rho,\sigma
\rrbracket\backslash q([0,1])$, and in particular, $\eta\not
=\rho,\sigma$. Then
$q([0,1])$ is contained in the union of the two disjoint open
sets $\t_g\backslash \t_g[\eta]$
and $\t_g[\eta]\backslash\{\eta\}$, with 
$q(0)=\rho\in\t_g\backslash\t_g[\eta]$ and $q(1)=\sigma
\in\t_g[\eta]\backslash\{\eta\}$.
This contradicts the fact that $q([0,1])$ is connected.

Conversely, suppose that there exists $a\in(0,1)$ such
that $q(a)\notin \llbracket\rho,\sigma\rrbracket$.
Set $\eta=q(a)$ and let $\gamma=\sigma\wedge \eta$. Note that $\gamma\in\llbracket \rho,\eta
\rrbracket \cap \llbracket \eta,\sigma
\rrbracket$ (from the definition
of $\sigma\wedge\eta$, it is immediate to verify that
$d_g(\eta,\sigma)=d_g(\eta,\gamma)+d_g(\gamma,\sigma)$). From the
first part of the proof of (ii),
$\gamma\in q([0,a])$ and, via a root change argument, $\gamma\in
q([a,1])$. Since $q$ is
injective, this is only possible if $\gamma=q(a)=\eta$, which contradicts
the
fact that $\eta\notin \llbracket\rho,\sigma\rrbracket$.
\cq

\medskip
\medskip
Once we know that $(\t_g,d_g)$ is an $\R$-tree, it is straightforward 
to verify that the notation $\sigma\preccurlyeq \sigma'$, $\llbracket \sigma,
\sigma'\rrbracket$, $\sigma\wedge \sigma'$ introduced in the preceding
proof is consistent with the definitions stated for a
general $\R$-tree at the beginning of this section.

\smallskip
Let us briefly discuss multiplicities of vertices in the tree
$\t_g$.
If $\sigma\in\t_g$ is not a leaf then we must have $\ell(\sigma)<r(\sigma)$,
where
$$\ell(\sigma):=\sup\,p_g^{-1}(\{\sigma\})\ ,\quad
r(\sigma):=\inf\,p_g^{-1}(\{\sigma\})$$
are respectively the smallest and the largest element in the
equivalence class
of
$\sigma$ in $[0,\zeta]$. Note that $m_g(\ell(\sigma),r(\sigma))
=g(\ell(\sigma))=g(r(\sigma))=d_g(\rho,\sigma)$.
Denote by $(a_i,b_i),i\in \ii$ the connected components of the open
set
$(\ell(\sigma),r(\sigma))
\cap\{t\in[0,\infty):g(t)>d_g(\rho,\sigma)\}$ (the index set $\ii$
is empty if $\sigma$ is a leaf). Then we claim that the connected
components  of the open set $\t_g\backslash\{\sigma\}$ are the sets
$p_g((a_i,b_i))$, $i\in \ii$ and $\t_g\backslash \t_g[\sigma]$ (the latter
only if $\sigma$ is not the root).  We have already noticed
that $\t_g\backslash \t_g[\sigma]$ is open, and the argument
used above for $\t_g[\sigma]\backslash\{\sigma\}$ also shows that the sets
$p_g((a_i,b_i))$,
$i\in \ii$ are open. Finally the sets $p_g((a_i,b_i))$ are
connected as continuous images of intervals, and
$\t_g\backslash \t_g[\sigma]$ is also connected because
if $\sigma',\sigma''\in \t_g\backslash \t_g[\sigma]$,
$\llbracket\rho,\sigma'\rrbracket
\cup \llbracket\rho,\sigma''\rrbracket$
is a connected closed set contained in $\t_g\backslash \t_g[\sigma]$.

\subsection{Convergence of trees}

Two rooted $\R$-trees $\t_{(1)}$ and $\t_{(2)}$ are called equivalent if
there is a root-preserving isometry that maps 
$\t_{(1)}$ onto $\t_{(2)}$. We denote by ${\T}$ the set of all 
equivalence classes of rooted compact $\R$-trees. The set $\T$
can be equipped with the (pointed) Gromov-Hausdorff distance, which
is defined as follows. 

If $(E,\delta)$ is a metric space, we use the
notation $\delta_{Haus}(K,K')$ for the usual Hausdorff metric between
compact subsets of $E$.
Then, if $\t$ and $\t'$ are two rooted compact $\R$-trees
with respective roots $\rho$ and $\rho'$, we define the distance
$d_{GH}(\t,\t')$ as 
$$d_{GH}(\t,\t')=\inf\Big(\delta_{Haus}(\varphi(\t),\varphi'(\t'))\vee \delta(\varphi(\rho),
\varphi'(\rho'))\Big),$$ 
where the infimum is over all isometric embeddings $\varphi:\t\la E$ and
$\varphi':\t'\la E$ of $\t$ and $\t'$ into a common metric
space $(E,\delta)$. Obviously $d_{GH}(\t,\t')$ only depends on
the equivalence classes of $\t$ and $\t'$. Furthermore $d_{GH}$ defines 
a metric on $\T$ (cf \cite{Gro} and \cite{EPW}).

According to Theorem 2 of \cite{EPW}, the metric space $({\T},d_{GH})$
is complete and separable. Furthermore, the distance 
$d_{GH}$ can often be evaluated in the following way. First recall that if
$(E_1,d_1)$ and $(E_2,d_2)$ are two compact metric spaces, a
correspondence between $E_1$ and $E_2$ is a subset ${\cal R}$ of
$E_1\times E_2$ such that for every $x_1\in E_1$ there exists at least one
$x_2\in E_2$ such that $(x_1,x_2)\in{\cal R}$ and conversely
for every $y_2\in E_2$ there exists at least one
$y_1\in E_1$ such that $(y_1,y_2)\in{\cal R}$. The distorsion of
the correspondence ${\cal R}$ is defined by
$${\rm dis}({\cal R})=\sup\{|d_1(x_1,y_1)-d_2(x_2,y_2)|:
(x_1,x_2),(y_1,y_2)\in{\cal R}\}.$$
Then, if $\t$ and $\t'$ are two rooted $\R$-trees with respective roots
$\rho$ and $\rho'$, we have
\begin{equation}
\label{bounddist}
d_{GH}(\t,\t')={1\over 2}\ \inf_{{\cal R}\in{\cal
C}(\t,\t'),\,(\rho,\rho')\in{\cal R}}\,{\rm dis}({\cal R}),
\end{equation}
where ${\cal
C}(\t,\t')$ denotes the set of all correspondences 
between $\t$ and $\t'$ (see Lemma 2.3 in \cite{EPW}).

\begin{lemma}
\label{dist-trees}
Let $g$ and $g'$ be two continuous functions with compact support
from $[0,\infty)$ into $[0,\infty)$, such that $g(0)=g'(0)=0$.
Then,
$$d_{GH}(\t_g,\t_{g'})\leq 2\|g-g'\|,$$
where $\|g-g'\|$ stands for the uniform norm of $g-g'$.
\end{lemma}

\noindent{\bf Proof.} We can construct a correspondence between 
$\t_g$ and $\t_{g'}$ by setting
$${\cal R}=\{(\sigma,\sigma'):\sigma=p_g(t)\hbox{ and }
\sigma'=p_{g'}(t)\hbox{ for some }t\geq 0\}.$$
In order to bound the distortion of $\cal R$, let
$(\sigma,\sigma')\in{\cal R}$ and $(\eta,\eta')\in{\cal R}$. By
our definition of ${\cal R}$ we can find $s,t\geq 0$ such that
$p_g(s)=\sigma$, $p_{g'}(s)=\sigma'$ and $p_g(t)=\eta$, $p_{g'}(t)=
\eta'$. Now recall that
\begin{eqnarray*}
d_g(\sigma,\eta)&=&g(s)+g(t)-2m_g(s,t),\\
d_{g'}(\sigma',\eta')&=&g'(s)+g'(t)-2m_{g'}(s,t),
\end{eqnarray*}
so that
$$|d_g(\sigma,\eta)-d_{g'}(\sigma',\eta')|\leq 4\|g-g'\|.$$
Thus we have ${\rm dis}({\cal R})\leq 4\|g-g'\|$ and the desired result
follows from (\ref{bounddist}). \cq

\section{The height process}

\subsection{The definition of the height process}

We will now introduce the random process which codes, in the sense of
subsection 2.1, the genealogical structure of a continuous-state branching process.
Recall that a continuous-state branching process is a Markov process $(Y_t,t\geq 0)$
with values in the positive half-line $[0,\infty)$, with a Feller semigroup
$(Q_t,t\geq 0)$ satisfying the following additivity (or branching) property: For every $t\geq 0$
and $x,x'\geq 0$,
$$Q_t(x,\cdot)*Q_t(x',\cdot)=Q_t(x+x',\cdot).$$
Informally, this is just saying that the union of two independent populations 
started respectively at $x$ and $x'$ will evolve like a single population started
at $x+x'$.

We will consider only the critical or subcritical case, meaning that
$\int_{[0,\infty)} y\,Q_t(x,dy)\leq x$ for every $t\geq 0$ and $x\geq 0$. Then the 
Laplace functional of the semigroup can be written in the following form:
\begin{equation}
\label{LaplaceCSBP}
\int_{[0,\infty)} e^{-\lambda y}\,Q_t(x,dy)=\exp(- x\,u_t(\lambda)),
\end{equation}
where the function $(u_t(\lambda),t\geq 0,\lambda\geq 0)$ is determined 
by the differential equation
\begin{equation}
\label{eqCSBP}
{du_t(\lambda)\over dt}=-\psi(u_t(\lambda))\ ,\quad u_0(\lambda)=\lambda\;,
\end{equation}
and $\psi$ is a function of the type
$$ \psi (\lambda )= \alpha \lambda + \beta \lambda^2 + 
\int_{(0, \infty)} (e^{-\lambda r} -1+\lambda r)\,\pi (dr)  \; ,  $$
where $\alpha, \beta \geq 0$ and $\pi $ is a $\sigma $-finite measure on $(0, \infty)$ such that 
$\int_{(0, \infty)} (r\wedge r^2)\,\pi (dr) < \infty $. Conversely, for any function $\psi$
of this type, there exists a (unique in law) continuous-state branching process $Y$
whose transition kernel is determined from $\psi$ by the preceding formulas. The process
$Y$ is called the $\psi$-continuous-state branching process ($\psi$-CSBP in short). It is well
known that $Y$ has only positive jumps (indeed $Y$ can be obtained as a time change of
a	spectrally positive L\'evy process, see Lamperti \cite{Lam}).

In the present work, we will consider only the case where 
the $\psi$-CSBP becomes extinct almost surely, which is equivalent to the condition
\begin{equation}
\label{extinct}
\int_1^\infty {du\over \psi(u)}<\infty.
\end{equation}
Note that this implies that at least one of the following two conditions holds: 
\begin{equation}
\label{infvar}
\beta >0 \quad {\rm or} \quad \int_{(0, 1)} r\,\pi (dr)  = \infty .
\end{equation}
(\ref{infvar}) is necessary and sufficient for the paths of $Y$
to be of infinite variation a.s. 
The coding of the genealogy that is presented below remains valid under (\ref{infvar})
even if (\ref{extinct}) fails to hold, but the resulting
tree is no longer compact (see Theorem 4.7 in \cite{LGLJ1}). 
On the other hand, if (\ref{infvar}) is not satisfied 
(that is in the finite variation case), the underlying
branching structure is basically discrete: See Section 3 of \cite{LGLJ1} and also
\cite{Li} for  a discussion 
with applications to queuing processes). 

Special cases that satisfy our assumptions are the quadratic
branching case $\psi(u)=c\,u^2$ and the stable branching case $\psi(u)=c\,u^{\gamma}$, for some
$1<\gamma<2$.

\medskip
It has been argued in \cite{LGLJ1} and \cite{DuLG} that the genealogy of the $\psi$-CSBP is coded
by the so-called height process, which is itself a functional of the 
L\'evy process with Laplace exponent $\psi$. We denote by $X$ a (spectrally positive) 
L\'evy process with Laplace
exponent $\psi$, defined under the probability measure $\P$:
$$ \E [ \exp ( -\lambda X_t ) ]= \exp (t\psi (\lambda )) \; ,\quad t, \lambda \geq 0 \; .$$
The subcriticality assumption on $Y$ and condition (\ref{infvar})
are equivalent to saying respectively that $X$ does no drift to $+\infty$
and has paths of infinite variation.

The height process $H=(H_t ; \; t\geq 0)$ associated with $X$ is
defined in such a way that, for every $t\geq 0$, $H_t$ measures the size of the set
\begin{equation}
\label{set}
\{ s\in [0,t] \; :\;  X_{s-} \leq\inf_{s\leq r \leq t} X_r \} \; .
\end{equation}
This is motivated by a discrete analogue for Galton-Watson trees (see Section 0.2 in \cite{DuLG}).
To make the preceding definition precise, we use a time-reversal argument: For any $t>0$, we define the
L\'evy process reversed at time $t$ by 
$$ \widehat{X}^t_s = X_t -X_{(t-s)-} \; ,\quad 0\leq s <t \quad {\rm and } \quad  \widehat{X}^t_t =X_t .$$
Then $\widehat{X}^t$ is distributed as $X$ up to time $t$. Let us set 
$$ S_s = \sup_{r\leq s} X_r \quad {\rm and } \quad \widehat{S}^t_s = \sup_{r\leq s} \widehat{X}^t_r .$$
The set (\ref{set}) is the image of
$$ \{ s\in [0,t] \; :\;  \widehat{S}^t_s =\widehat{X}^t_s  \} \; .$$
under the time reversal operation $s\to t-s$. Recall that 
$S-X$ is a strong Markov process for which $0$ is a regular point. So, we can consider its local time process at $0$,
which is denoted by 
$\Gamma(X)=(\Gamma_t(X),t\geq 0)$. We define the height process by 
\begin{equation}
\label{defheight}
H_t = \Gamma_t(\widehat{X}^t ) \; ,\quad t\geq 0 .
\end{equation}

To complete the definition, we still need to specify the normalization of the local time
$\Gamma(X)$. This can be done
through the
following approximation:
$$H_t=\lim_{\varepsilon\da 0} {1\over \varepsilon}\int_0^t ds\,\un_{\{X_s\leq I^s_t+\varepsilon\}},$$
where $I^s_t=\inf_{s\leq r\leq t}X_r$ 
and the convergence holds in probability (this approximation follows from Lemma 1.1.3 in \cite{DuLG}).
Thanks to  condition (\ref{extinct}), we know that the process $H$ has a modification with 
continuous sample paths (Theorem 4.7 in \cite{DuLG}). From now on we consider only this 
modification. When
$\beta >0$, it is not hard to see that, for any $t\geq 0$,
$$H_t= \frac{1}{\beta} {\rm Leb} \left( \{ \widehat{S}^t_s \; ;\; 0\leq s \leq t\} \right) , $$
where $ {\rm Leb}$ stands for the Lebesgue measure on the real line.
In particular when $\psi(u)=\beta u^2$ ($X$ is a scaled Brownian motion), we see that
$H_t=\beta^{-1}(S_t-X_t)$ is distributed as a  (scaled) reflected
Brownian motion.

\medskip
For our purposes it will be crucial to define the height process also under the
so-called excursion measure $N$.
Set $I_t=\inf_{s\leq t} X_s $ and recall that $X-I$ is a strong Markov process.
Then for any $t\geq 0$, $H_t$ only depends on the values taken by $X-I$
on the excursion interval that straddles $t$ (at least informally this is obvious if we think
of $H_t$ as measuring the size of the set (\ref{set})). Under our assumptions, $0$ is
a regular point for $X-I$, and the process $-I$ can be chosen as the local time of $X-I$ at level $0$. 
We denote by $N$
the associated excursion measure, which plays a fundamental role throughout
this work (as was already the case in \cite{DuLG}). The duration of the excursion under $N$
is denoted by $\zeta$. Let $(g_i,
d_i)$,
$i\in\ii$ be the excursion intervals of $X-I$ above $0$. One easily verifies that $\P$ a.s.,
$$ \bigcup_{i\in \ii} (g_i , d_i) = \{ s\geq 0 \; :\;  X_s -I_s >0 \} =\{  s\geq 0 \; :\;  H_s >0 \} .$$
Denote by $H_i (s) =H_{g_i+s}$ , $ 0\leq s\leq \zeta_i=d_i-g_i $, $i\in \ii$ the excursions 
of $H$ away from $0$. Then, each
$H_i$ can be written as a functional of the excursion of 
$X-I$ away from $0$ corresponding to the interval $(g_i,d_i)$.
Consequently, if $C_+([0,\infty))$ denotes the space of all 
nonnegative continuous functions on $[0,\infty)$, the point measure 
\begin{equation}
\label{Poiss}
\sum_{i\in \ii } \delta_{(-I_{g_i}, H_i )} (d\ell d\omega)
\end{equation}
is a Poisson point measure on $\R_+\times C_+([0,\infty))$ with intensity $d\ell\Delta(d\omega)$, where
$\Delta(d\omega)$  is the $\sigma$-finite measure on $C_+([0,\infty))$ defined as the law
of $H$ under
$N$. Note that in the Brownian case, $\Delta$ is the classical It\^o measure of positive excursions
of linear Brownian motion (up to a normalizing constant).

\subsection{Local times of the height process}

Let us start by the defining the local times under $\P$. For every 
$a\geq 0$, the local time of the height process at level $a$ is the
continuous increasing process $(L^a_s,s\geq 0)$ which can be characterized via the
approximation
$$ \lim_{\varepsilon \rightarrow 0} \E \left[  \sup_{ 0\leq s \leq t} \left| \frac{1}{\varepsilon} \int_0^s dr 
\un_{\{ a< H_r \leq a+\varepsilon \}} -L_s^a \right| \right] =0 $$
(see Section 1.3 of \cite{DuLG}). It is then easy to see that
the support of the measure $dL^a_s$ is contained in the closed set $\{s\geq 0:H_s=a\}$. When $a>0$, we have
also
$$ \lim_{\varepsilon \rightarrow 0} \E \left[  \sup_{ 0\leq s \leq t} \left| \frac{1}{\varepsilon} \int_0^s dr 
\un_{\{ a-\varepsilon< H_r \leq a\}} -L_s^a \right| \right] =0\;. $$

 Let us recall the ``Ray-Knight theorem'' for $H$ 
(\cite{DuLG} Theorem 1.4.1, see also \cite{LGLJ1}, Theorem 4.2),
which can be viewed as a generalization of famous 
results about linear Brownian motion. For any $r\geq 0$, set: $T_r=\inf \{ s\geq 0 \; :\;
X_s=-r\}$. Then, the process $(L^a_{T_r} \; ;\; a\geq 0)$ is a $\psi$-CSBP started at $r$.
In particular, this process has a c\` adl\` ag modification.

The local time at level $a$ can also be used to describe the distribution of 
excursions of the height process above level $a$, and this will be very important 
for our applications. Let us fix $a>0$ and denote by $(\alpha_j,\beta_j)$, $j\in \jj$
the connected components of the open set $\{s\geq 0:H_s>a\}$. For any $j\in \jj$,
denote by $H^j$ the corresponding excursion of $H$ defined by:
$$H^j_s=H_{(\alpha_j+s)\wedge \beta_j}-a\ ,\quad s\geq 0.$$
Also set $\wt H^a_s=H_{\tilde\tau^a_s}$, where for 
every $s\geq 0$,
$$\wt\tau^a_s=\inf\{t\geq 0:\int_0^t dr\,\un_{\{H_r\leq a\}}>s\}.$$
Informally, $\wt H^a$ corresponds to the evolution of $H$ ``below level $a$''.

The next result is a straightforward consequence of Proposition 1.3.1 in \cite{DuLG}.

\begin{proposition}
\label{subexclevel}
Under the probability ${\bf P}$, the point measure 
$$\sum_{j\in \jj} \delta_{(L^a_{\alpha_j},H^j)}(d\ell \,d\omega)$$
is independent of $\wt H^a$ and is a Poisson point measure on $\R_+\times C_+([0,\infty))$ with intensity 
$d\ell \Delta(d\omega)$.
\end{proposition}

It will be important to define local times under the excursion measure $N$. This creates no
additional difficulty thanks to the following simple remark. If $r>0$, then for any $\delta>0$, there is a
positive probability under $\P$ that
exactly one excursion of $H$ away from zero hits level $\delta$ before time $T_r$. It easily follows that
we can define for every $a>0$ a continuous increasing process $(\Lambda^a_s,s\geq 0)$, such that, for every
$\delta\in(0,a)$ and $t\geq 0$,
\begin{equation}
\label{localapprox}
\lim_{\varepsilon \rightarrow 0} \, 
N \left(  \un_{\{\sup H>\delta\}}\;\sup_{ 0\leq s \leq t\wedge \zeta} \left| \frac{1}{\varepsilon} \int_0^s
dr 
\un_{\{ a-\varepsilon< H_r \leq a\}} -\Lambda_s^a \right| \right) =0.
\end{equation}
(see Section 1.3 in \cite{DuLG}). Again the support of the measure 
$d\Lambda^a_s$ is contained in $\{s:H_s=a\}$, $N$ a.e. From the
above-mentioned
Ray-Knight  theorem and elementary excursion theory for $X-I$ we get, for any $a>0$ and any $
\lambda\geq 0$, 
\begin{equation}
\label{localloi}
N\left( 1-\exp (-\lambda \Lambda^a_{\zeta} ) \right) =u_a (\lambda ),
\end{equation}
where $u_a(\lambda)$ is as in (\ref{LaplaceCSBP}).
We set $ v(a) = \lim_{\lambda \rightarrow \infty} 
u_a (\lambda)$. By writing (\ref{eqCSBP}) in the form of an integral equation
and passing to the limit $\lambda\to\infty$ using (\ref{extinct}), we see that the 
function $v$ is finite on $(0,\infty)$ and determined by the equation
$$a= \int_{v(a)}^{\infty} {du\over\psi (u)}.$$ 
Moreover, for every $a>0$, we have
\begin{equation}
\label{totheig}
v(a)=N\left(  \Lambda^a_{\zeta}>0 \right) =N\left(  \sup_{s\leq \zeta} H_s >a \right). 
\end{equation} 
The first equality follows from the definition of $v$. The second one can
be deduced from Proposition \ref{subexclevel}, which implies that
$\inf\{s\geq 0:L^a_s>0\}=\inf\{s\geq 0:H_s>a\}$, $\P$ a.s.

\smallskip
We will need an analogue of Proposition 3.1 under the excursion measure $N$. To state
it, fix $a>0$ and denote by $(\alpha_j,\beta_j)$, $j\in \jj$ the excursion intervals of $H$ above level
$a$ (just as before, but we are now arguing under $N$) and for every $j\in \jj$
let $H^j$ be the corresponding excursion. Let the process $\wt H^a$
be defined as previously and let $\wt\h^a$ be the $\sigma$-field generated
by $\wt H^a$ and the class of $N$-negligible measurable sets. From our approximation 
(\ref{localapprox}) it follows that $\Lambda^a_\zeta$ is measurable with respect
to $\wt \h^a$.

\begin{corollary}
\label{sublevel}
Under the probability measure $N(\cdot\mid \sup H>a)$ and conditionally on $\wt \h^a$, the point measure
$$\sum_{j\in \jj} \delta_{(\Lambda^a_{\alpha_j},H^j)}(d\ell \,d\omega)$$
is distributed as a Poisson point measure on $\R_+\times C_+([0,\infty))$ with intensity 
$\un_{[0,\Lambda^a_\zeta]}(\ell)d\ell \Delta(d\omega)$.
\end{corollary}

This is really an immediate consequence of Proposition \ref{subexclevel} if we notice that the
law under $\P$ of the
first excursion of $H$ that hits level $a$ is $N(\cdot\mid \sup H>a)$. Alternatively, the statement
of Corollary \ref{sublevel} also appears as an intermediate result in the proof of Proposition 4.2.3 in
\cite{DuLG}.

We will need one additional
property related to Corollary \ref{sublevel}. First denote by $(\wt\Lambda^a_s,s\geq 0)$
the local time of $\wt H^a$ at level $a$, which may be defined either by
an approximation similar to (\ref{localapprox}) or directly by the formula 
$\wt\Lambda^a_s=\Lambda^a_{\tilde
\tau^a_s}$. Then we have $N$ a.e. on $\{\sup H>a\}$
\begin{equation}
\label{techsub}
\inf\{s\geq 0:\wt \Lambda^a_s >\Lambda^a_{\alpha_j}\}=\int_0^{\alpha_j} ds\,\un_{\{H_s\leq a\}}\ ,
\hbox{ for every }j\in \jj.
\end{equation}
For a proof, see pages 108-109 of \cite{DuLG}.

We conclude this section with an important regularity property of local times. Recall that a
c\` adl\` ag process $Y$ is said to have no fixed discontinuities if for every fixed $t>0$,
the sample path of $Y$ is continuous at $t$ outside a set of zero probability.

\begin{lemma}
\label{contiLT}
Set $\Lambda^0_s=0$ for every $s\geq 0$. Then the process $(\Lambda^a_\zeta,a\geq 0)$
has a c\` adl\` ag modification under $N$, and this modification has no fixed discontinuities.
\end{lemma}

\noindent{\bf Proof.} Let $r>0$. Since the process $(L^a_{T_r},a\geq 0)$ is 
a $\psi$-CSBP and thus a Feller process, it has a c\` adl\` ag modification 
with no fixed discontinuities under $\P$. Let $H_i$, $i\in \ii$
be the excursions of $H$ away from $0$, as in (\ref{Poiss}), and for
every $i\in \ii$ let $\zeta_i$ be the duration of $H_i$. From our
approximation of local times, it is easy to see that, for every $a>0$,
\begin{equation}
\label{localtimeDecomp}
L^a_{T_r}=\sum_{i\in\ii,I_{g_i}>-r} \Lambda^a_{\zeta_i}(H_i)\ ,\quad N\hbox{ a.e.}
\end{equation}
Using a previous remark about the existence of exactly one excursion of $H$ hitting level
$\delta$ before time $T_r$, we easily deduce from the previous
formula and the c\` adl\` ag property of $(L^a_{T_r},a\geq 0)$ that the process
$(\Lambda^a_{\zeta},a>0)$ must have a c\` adl\` ag modification with no
fixed discontinuities under $N$. 
Furthermore, if we use this modification in the right side of (\ref{localtimeDecomp}), for 
every $a>0$,
we will obviously obtain the c\` adl\` ag modification of the process $(L^a_{T_r},a> 0)$.

It remains 
to verify that $\Lambda^a_\zeta$ converges to $0$,
$N$ a.e. as $a\da 0$ (we now consider the modification that has just been introduced). For this, we need a
different argument. Let
$\delta>0$ and let
$H_{i_0}$ be the first excursion of $H$ that reaches level $\delta$. From 
properties of Poisson measures, the law under $\P$ of the point measure
$$\sum_{i\in \ii\backslash\{i_0\},I_{g_i}>-r } \delta_{(-I_{g_i}, H_i )} (drd\omega)$$
is absolutely continuous with respect to that of
$$\sum_{i\in \ii,I_{g_i}>-r  } \delta_{(-I_{g_i}, H_i )} (drd\omega).$$
In particular,
the function 
$$a\la \sum_{i\in\ii\backslash\{i_0\},I_{g_i}>-r} \Lambda^a_{\zeta_i}(H_i)$$
must converge $\P$ a.s. to $r$ as $a\da 0$. Now note that, on the event $\{-I_{g_{i_0}}>-r\}
=\{\sup_{[0,T_r]}H>\delta\}$, we have for every $a>0$
$$\Lambda^a_{\zeta_{i_0}}(H_{i_0})=
\sum_{i\in\ii,I_{g_i}>-r} \Lambda^a_{\zeta_i}(H_i)
-\sum_{i\in\ii\backslash\{i_0\},I_{g_i}>-r} \Lambda^a_{\zeta_i}(H_i)$$
and use the fact that the distribution of $H_{i_0}$ under $\P(\cdot\mid \sup_{[0,T_r]}H>\delta)$
coincides with the law of $H$ under $N(\cdot \mid \sup H>\delta)$ to complete the proof. \cq

\smallskip
From now on, we assume that have chosen a modification of the collection
$(\Lambda^a,a\geq 0)$ in such a way that the process $(\Lambda^a_\zeta,a\geq 0)$
is c\` adl\` ag. This will be important in the applications developed in
Section 4 below.
  
Let us finally briefly comment on the use of the measures $\P$ and $N$
for our purposes. As will be made precise in the next section, the height process
under $N$ codes a single (compact rooted) $\R$-tree, whereas under $\P$ it codes
a Poissonnian collection of such trees, each excursion of $H$ away from $0$
corresponding to one tree.

\subsection{A key lemma}

In this subsection, we prove a basic preliminary lemma, which is a consequence of
the results in \cite{DuLG}. We need to introduce some notation. Denote by
$M_f$ the space of all finite measures on $[0,\infty)$. If $\mu\in M_f$, we 
denote by $H(\mu)\in[0,\infty]$ the
supremum  of the (topological) support of $\mu$. We also introduce a ``killing operator''
on measures defined as follows. For every $x\geq 0$, $k_x\mu$ is the element of $M_f$ such that
$k_x\mu([0,t])=\mu([0,t])\wedge (\mu([0,\infty))-x)_+$ for every $t\geq 0$. Let $M^*_f$ stand for
the set of all measures  $\mu\in M_f$
such that $H(\mu)<\infty$ and the topological support of $\mu$
is $[0,H(\mu)]$. If $\mu\in M^*_f$, we denote by $Q_\mu$ the law under $\P$ of the process $H^\mu$
defined by
$$
\begin{array}{ll}
H^\mu_t=H(k_{-I_t}\mu)+H_t\;,&\quad \hbox{if }t\leq T_{\langle \mu,1\rangle}\;,\\
H^\mu_t=0\;,&\quad \hbox{if }t>T_{\langle \mu,1\rangle}\;,
\end{array}
$$
where $T_{\langle \mu,1\rangle}=\inf\{t\geq 0:X_t=-\langle \mu,1 \rangle\}$.
Our assumption $\mu\in M^*_f$ guarantees that $H^\mu$ has continuous sample paths, and
we can therefore view $Q_\mu$ as a probability measure on the space 
$C_+([0,\infty))$ of nonnegative continuous functions on $[0,\infty)$. 

Finally, let $\psi^*(u)=\psi(u)-\alpha\,u$, and let 
$(U^1,U^2)$ be a two-dimensional subordinator with Laplace functional
$$E[\exp(-\lambda U^1_t-\lambda'U^2_t)]=\exp\Big(-{\psi^*(\lambda)-\psi^*(\lambda')\over \lambda-\lambda'}
\Big).$$
(alternatively, $(U^1,U^2)$ can be characterized by its drift and L\'evy measure,
see \cite{DuLG}, p. 80). When $\lambda=\lambda'$, the ratio ${\psi^*(\lambda)-\psi^*(\lambda')\over
\lambda-\lambda'}$ should obviously be interpreted as 
$\psi'(\lambda)-\alpha$, so that we see that $U^1+U^2$ is
a subordinator with Laplace exponent $\psi'-\alpha$. For every $a\geq 0$, we let $\M_a$ be the
probability measure on
$(M^*_f)^2$ which is the distribution of $(\un_{[0,a]}(t)\,dU^1_t,\un_{[0,a]}(t)\,dU^2_t)$.

\begin{lemma}
\label{keytech}
For any nonnegative measurable function $F$ on $C_+([0,\infty))^2$,
\begin{eqnarray*}
&&N\Big(\int_0^\zeta ds\,F\Big((H_{(s-t)_+},t\geq 0),(H_{(s+t)\wedge \zeta},t\geq 0)\Big)\Big)\\
&&\quad=\int_0^\infty da\,e^{-\alpha a}\int \M_a(d\mu d\nu)\int Q_\mu(dh)Q_\nu(dh')F(h,h').
\end{eqnarray*}
\end{lemma}

\noindent{\bf Remark.} In the Brownian case $\psi(u)=u^2$, this lemma reduces to the 
well-known Bismut decomposition of the Brownian excursion.

\smallskip

\noindent{\bf Proof.} We start by recalling some results from \cite{DuLG} (see 
Chapter 1 and Section 3.1 in \cite{DuLG}). We can define both under $\P$ 
and under $N$ a c\` adl\` ag process $(\rho_t,\eta_t)_{t\geq 0}$ taking values in $(M_f^*)^2$
such that the following properties hold:
\begin{description}
\item{(i)} We have $H(\rho_t)=H_t=H(\eta_t)$ for every $t\geq 0$, $N$ a.e. and $\bf P$ a.e.
\item{(ii)} The process $(\rho_t,\eta_t)$ is adapted with respect to the filtration
$({\cal F}_t)_{t\geq 0}$ generated
by the L\'evy process $X$. Furthermore, if $G$ is any nonnegative measurable functional 
on $C_+([0,\infty))$, we have for every $s>0$, $N$ a.e. on the event $\{s<\zeta\}$,
\begin{equation}
\label{potker}
N\Big(G(H(\rho_{(s+t)\wedge \zeta}),t\geq 0)\;\Big|\;{\cal F}_s\Big)=Q_{\rho_s}(G).
\end{equation}
\item{(iii)} The process $(\eta_{(\zeta-s)_-},\rho_{(\zeta-s)_-})_{0\leq s<\zeta}$
has the same distribution as $(\rho_s,\eta_s)_{0\leq s<\zeta}$ under $N$.
\item{(iv)} For any nonnegative measurable function $\Phi$ on $(M_f^*)^2$,
\begin{equation}
\label{invar}
N\Big(\int_0^\zeta ds\,\Phi(\rho_s,\eta_s)\Big)=\int_0^\infty da\,e^{-\alpha a}\,\M_a(G).
\end{equation}
\end{description}
To make sense of the conditional expectation in (\ref{potker}), note that the event
$\{s<\zeta\}$ has finite $N$-measure. We refer to \cite{DuLG} for a proof of properties
(i) -- (iv) above: See in particular Propositions 1.2.3 and 1.2.6 for (\ref{potker}),
Corollary 3.1.6 for (iii) and Proposition 3.1.3 for (iv).

We now proceed to the proof of the lemma. We may and will assume that 
$F$ is of the form $F(h,h')=F_1(h)F_2(h')$. Using (i) and then (ii), we have
\begin{eqnarray*}
&&N\Big(\int_0^\zeta ds\,F_1\Big(H_{(s-t)_+},t\geq 0\Big)F_2\Big(H_{(s+t)\wedge \zeta},t\geq 0\Big)\Big)\\
&&\quad=N\Big(\int_0^\zeta ds\,F_1\Big(H(\rho_{(s-t)_+}),t\geq 0\Big)
F_2\Big(H(\rho_{(s+t)\wedge \zeta}),t\geq
0\Big)\Big)\\
&&\quad=N\Big(\int_0^\zeta ds\,F_1\Big(H(\rho_{(s-t)_+}),t\geq 0\Big)\,Q_{\rho_s}(F_2)\Big).
\end{eqnarray*}
From the time-reversal property (iii) we see that the last quantity is equal to
\begin{eqnarray*}
&&N\Big(\int_0^\zeta ds\,Q_{\eta_s}(F_2)\,F_1\Big(H(\rho_{(s+t)\wedge \zeta}),t\geq
0\Big)\Big)\\
&&\quad=N\Big(\int_0^\zeta ds\,Q_{\eta_s}(F_2)\,Q_{\rho_s}(F_1)\Big),
\end{eqnarray*}
using (\ref{potker}) once again. The formula of the lemma now follows from (\ref{invar}). \cq

\begin{corollary} 
\label{localkey}
Let $a>0$. Then, for any nonnegative measurable function $F$ on $C_+([0,\infty))^2$,
\begin{eqnarray*}
&&N\Big(\int_0^\zeta d\Lambda^a_s\,F\Big((H_{(s-t)_+},t\geq 0),(H_{(s+t)\wedge \zeta},t\geq 0)\Big)\Big)\\
&&\quad=e^{-\alpha a}\int \M_a(d\mu d\nu)\int Q_\mu(dh)Q_\nu(dh')F(h,h').
\end{eqnarray*}
\end{corollary}

\noindent{\bf Proof.} 
This is a straightforward consequence of Lemma \ref{keytech} and the approximation of local time
given in (\ref{localapprox}). \cq

\medskip
\noindent{\bf Remark.} The case $F=1$ of Corollary \ref{localkey} gives $N(\Lambda^a_\zeta)=e^{-\alpha a}$,
for every $a>0$.

\section{The L\'evy tree}

We have seen that $N$ a.e. the function $s\to H_s$ satisfies the properties stated at
the beginning of subsection 2.1, namely it is continuous with compact support and
such that $H_0=0$. 

\begin{definition}
The $\psi$-L\'evy tree
is the tree $(\t_H,d_H)$
coded by the function $s\to H_s$ under the measure $N$.
\end{definition}

We will say the L\'evy tree rather than the $\psi$-L\'evy tree if there is no risk of confusion.

 We denote by $\Theta(d\t)$ the $\sigma$-finite 
measure on $\T$ which is the law 
of the L\'evy tree, that is the law of the tree $\t_H$ under $N$. Note that the 
measurability of the random variable  $\t_{H}$ follows from Lemma
\ref{dist-trees}.

\subsection{From discrete to continuous trees}

In this subsection, we will state a result which justifies the definition 
of the $\psi$-L\'evy tree by showing that it arises as the limit 
in the Gromov-Hausdorff distance of
suitably rescaled discrete Galton-Watson trees.

We start by introducing some formalism for discrete trees. Let
$${\cal U}=\bigcup_{n=0}^\infty \N^n  $$
where $\N=\{1,2,\ldots\}$ and by convention $\N^0=\{\varnothing\}$. If
$u=(u_1,\ldots u_m)$ and 
$v=(v_1,\ldots, v_n)$ belong to $\cal U$, we write $uv=(u_1,\ldots u_m,v_1,\ldots ,v_n)$
for the concatenation of $u$ and $v$. In particular $u\varnothing=\varnothing u=u$.

A (finite) rooted ordered tree $\theta$ is a finite subset of
$\cal U$ such that:
\begin{description}
\item{(i)} $\varnothing\in \theta$.

\item{(ii)} If $v\in \theta$ and $v=uj$ for some $u\in {\cal U}$ and
$j\in\N$, then $u\in\theta$.

\item{(iii)} For every $u\in\theta$, there exists a number $k_u(\theta)\geq 0$
such that $uj\in\theta$ if and only if $1\leq j\leq k_u(\theta)$.
\end{description}

\noi We denote by ${\bf T}$ the set of all rooted ordered trees. In what follows, we see each vertex of the
tree $\theta$ as an individual of a population  whose $\theta$ is the family tree. 

If $\theta$ is a tree and $u\in \theta$, we define the shift of $\theta$ at $u$
by $\tau_u\theta=\{v\in U:uv\in\theta\}$.
Note that $\tau_u\theta\in {\bf T}$. 
We also denote by $h(\theta)$ the height of $\t$, that is the maximal generation of a vertex in $\theta$.

\smallskip

\begin{center}
\unitlength=1.1pt
\begin{picture}(200,140)

\thicklines \put(50,0){\line(-1,2){20}}
\thicklines \put(50,0){\line(1,2){20}}
\thicklines \put(30,40){\line(0,1){40}}
\thicklines \put(30,40){\line(-1,2){20}}
\thicklines \put(30,40){\line(1,2){20}}
\thicklines \put(30,80){\line(-1,2){20}}
\thicklines \put(30,80){\line(1,2){20}}
\put(40,-5){$\varnothing$}
\put(20,35){$1$}
\put(60,35){$2$}
\put(-3,75){$11$}
\put(17,75){$12$}
\put(37,75){$13$}
\put(-7,115){$121$}
\put(31,115){$122$}
\put(30,-25){tree $\t^\theta$}
\thinlines \put(100,0){\line(1,0){120}}
\thinlines \put(100,0){\line(0,1){130}}
\thicklines \put(100,0){\line(1,5){8}}
\thicklines \put(108,40){\line(1,5){8}}
\thicklines \put(116,80){\line(1,-5){8}}
\thicklines \put(124,40){\line(1,5){8}}
\thicklines \put(132,80){\line(1,5){8}}
\thicklines \put(140,120){\line(1,-5){8}}
\thicklines \put(148,80){\line(1,5){8}}
\thicklines \put(156,120){\line(1,-5){8}}
\thicklines \put(164,80){\line(1,-5){8}}
\thicklines \put(172,40){\line(1,5){8}}
\thicklines \put(180,80){\line(1,-5){8}}
\thicklines \put(188,40){\line(1,-5){8}}
\thicklines \put(196,0){\line(1,5){8}}
\thicklines \put(204,40){\line(1,-5){8}}

\thinlines \put(108,0){\line(0,1){2}}
\thinlines \put(116,0){\line(0,1){2}}
\thinlines \put(124,0){\line(0,1){2}}
\thinlines \put(100,40){\line(1,0){2}}
\thinlines \put(100,80){\line(1,0){2}}
\thinlines \put(212,0){\line(0,1){2}}

\put(125,-25){contour function}

\put(107,-6){\scriptsize 1}
\put(115,-6){\scriptsize 2}
\put(123,-6){\scriptsize 3}

\put(95,39){\scriptsize 1}
\put(95,79){\scriptsize 2}

\end{picture}

\vskip 15mm

Figure 1
\end{center}

For our purposes it will be convenient to view $\theta$ as an $\R$-tree: To this end,
embed $\theta$ in the plane, in such a way that each edge corresponds to
a line segment of length one, in the way suggested by the left part of Fig. 1. Denote by $\t^\theta$ the
union of all these line segments  and equip $\t^\theta$ with the obvious metric such that the distance
between 
$\sigma$ and
$\sigma'$ is the length of the shortest path from $\sigma$ to $\sigma'$ in $\t^\theta$. This construction
leads to a (compact rooted) $\R$-tree whose
equivalence class does not depend on the particular embedding.

To define now the contour function of $\t^\theta$, consider a particle that starts from the root
and visits continuously  all edges at speed one, going backwards as less as possible and respecting 
the lexicographical order of vertices. Then let $C^\theta(t)$ denote the distance to the root
of the position of the particle at time $t$ (for $t\geq 2(|\theta|-1)$, we take $C^\theta(t)=0$
by convention).
Fig.1 explains
the definition of the contour function better than a formal definition. Note that in the notation of
Section 2, we have $\t^\theta=\t_{C^\theta}$, meaning that $\t^\theta$ coincides with the tree
coded by the function $C^\theta$.

Now let us turn to Galton-Watson trees.
Let $\mu$ be a critical or subcritical offspring distribution.  This means that
$\mu$ is a probability measure on $\Z_+$ such that
$\sum_{k=0}^\infty k\mu(k)\leq 1$.
We exclude the trivial case where $\mu(1)=1$.
Then,
there is a unique probability distribution $\Pi_\mu$ on ${\bf T}$ such that

\begin{description}
\item{(i)} $\Pi_\mu(k_\varnothing=j)=\mu(j)$,\quad
$j\in \Z_+$.

\item{(ii)} For every $j\geq 1$ with $\mu(j)>0$, the shifted trees
$\tau_1\theta,\ldots,\tau_j\theta$
are independent under the conditional probability
$\Pi_\mu(\cdot\mid k_\varnothing =j)$
and their conditional distribution is $\Pi_\mu$.
\end{description}

A random tree with distribution $\Pi_\mu$ is called a Galton-Watson tree
with offspring distribution $\mu$, or in short a $\mu$-Galton-Watson tree. Obviously
it describes the genealogy of the Galton-Watson process with offspring distribution $\mu$
started initially with one individual.

We can now state our result relating discrete and continuous trees. If $\t$
is a (compact rooted) $\R$-tree with metric $d$, and if $r>0$, we slightly abuse notation by writing
$r\t$ for the ``same'' tree equipped with the distance $r\,d$. Recall
that the height of $\t$ is 
$$h(\t)=\sup\{d(\rho(\t),\sigma):\sigma\in\t\},$$
where $\rho(\t)$ denotes the root of $\t$.  For every real number
$x$, $[x]$ denotes the integer part of $x$.

\begin{theorem}
\label{discrete-continuous}
Let $(\mu_p)_{p\geq 1}$ be  a sequence of critical or subcritical offspring distributions. 
For every $p\geq 1$ denote by $Y^p$ a Galton-Watson branching process with offspring
distribution $\mu_p$, started at $Y^p_0=p$. Assume that there exists a 
nondecreasing sequence $(m_p)_{p\geq 1}$ of positive integers converging to
$+\infty$ such that
\begin{equation}
\label{Grim1}
\left(p^{-1}Y^p_{[m_pt]}\,,\,t\geq 0\right)
\build{\longrightarrow}_{p\to\infty}^{{\rm (d)}} (Y_t,t\geq 0)
\end{equation}
where the limiting process $Y$ is a $\psi$-CSBP. Assume furthermore that
for every $\delta>0$,
\begin{equation}
\label{extinctBP}
\liminf_{p\to\infty} P[Y^p_{[m_p\delta]}=0]>0.
\end{equation}
Then, for every $a>0$, the law of the $\R$-tree $m_p^{-1}\t^\theta$ under $\Pi_{\mu_p}(d\theta\mid
h(\theta)\geq[a\,m_p])$ converges as $p\to\infty$ to the probability measure
$\Theta(d\t\mid h(\t)>a)$, in the sense of weak
convergence of measures in the space
$\T$.
\end{theorem}

\noindent{\bf Proof.} We noted that the tree
$\t^\theta$ is the tree 
coded by the function $C^\theta$, in the sense of Section 2. 
Lemma \ref{dist-trees} then shows that the convergence of
the theorem follows from the weak convergence of the scaled 
contour function $(m_p^{-1}C^\theta(pm_pt),t\geq 0)$ under $\Pi_{\mu_p}(d\theta\mid
h(\theta)\geq[a\,m_p])$ towards the height process $H$ under $N(\cdot\mid \sup H>a)$.
But this is precisely the contents of Proposition
2.5.2 in
\cite{DuLG}, which is itself a consequence of Theorem 2.3.1
in the same work. \cq

\medskip
The technical assumption (\ref{extinctBP}) guarantees that the
Galton-Watson process $Y^p$ dies out at a time of order $m_p$, as one expects
from the convergence (\ref{Grim1}) (recall that $Y$ dies out in finite time). See 
Chapter 2 of \cite{DuLG} for a discussion of this assumption, and note that it 
is always true in the case when $\mu_p=\mu$ for every $p$ (Theorem 2.3.2 in
\cite{DuLG}). In particular, the approximation result stated in the 
introduction above is easily seen to be a consequence of Theorem \ref{discrete-continuous}.

\subsection{Local times and the branching property of L\'evy trees}

Let us start with a few simple observations. We recall that the generic element of
$\T$ is denoted by $(\t,d)$. Then, for every $a>0$,
$$\Theta(h(\t)>a)=N(\sup H>a)=v(a),$$
where the function $v$ is determined by $\int_{v(a)}^\infty \psi(u)^{-1}du=a$.

The truncation of the tree $\t$ at level $a>0$ is the new tree
$${\rm tr}_a(\t)=\{\sigma\in \t:d(\rho(\t),\sigma)\leq a\},$$
which is obviously equipped with the restriction of the distance $d$. It is easy to
verify that the mapping $\t\to {\rm tr}_a(\t)$
from $\T$ into itself is measurable. 

Let $\t\in\T$ and $a>0$. Denote by $\t^{(i),\circ}$, $i\in \ii$ the connected 
components of the open set
$$\t({(a,\infty)})=\{\sigma\in \t:d(\rho(\t),\sigma) >a\}.$$
Notice that the index set $\ii$ may be empty (if $h(\t)\leq a$), finite or
countable. Let $i\in \ii$. Then the ancestor of $\sigma$ at level $a$ must be the
same for every $\sigma\in \t^{(i),\circ}$. We denote by $\sigma_i$ this 
common ancestor and set $\t^{(i)}=\t^{(i),\circ}\cup\{\sigma_i\}$. Then $\t^{(i)}$ is a compact rooted
$\R$-tree with root
$\sigma_i$. The trees $\t^{(i)}$, $i\in I$ are called the subtrees
of $\t$ originating from level $a$. 

We set
$$\n_a^\t:=\sum_{i\in \ii} \delta_{(\sigma_i,\t^{(i)})},$$
which is a point measure on $\t(a)\times \T$. Also, for every $\delta>0$, we let
$$Z(a,\delta):=|\{i\in \ii:h(\t^{(i)})\geq \delta\}|$$
be the number of subtrees of $\t$ originating from level $a$ that hit level $a+\delta$.

\begin{theorem}
\label{existLT}
For every $a\geq 0$ and for $\Theta$ a.e. $\t\in\T$
we can define a finite measure $\ell^a$ on $\t$, in such a way that the
following properties hold:
\begin{description}
\item{\rm (i)} $\ell^0=0$ and, for every $a>0$, $\ell^a$ is 
supported on $\t(a)$, $\Theta(d\t)$ a.e.
\item{\rm (ii)} For every $a>0$, $\{\ell^a\ne 0\}=\{h(\t)> a\}$, $\Theta(d\t)$ a.e.
\item{\rm (iii)} For every $a>0$, we have $\Theta(d\t)$ a.e. for every bounded
continuous function $\varphi$ on $\t$,
\begin{eqnarray}
\label{appLTlevel}
\langle\ell^a,\varphi\rangle&=&\lim_{\varepsilon\da 0}
{1\over v(\varepsilon)}\int \n^\t_a(d\sigma d\t')\,\varphi(\sigma)\,\un_{\{h(\t')\geq \varepsilon\}}
\nonumber\\
&=&\lim_{\varepsilon\da 0}
{1\over v(\varepsilon)}\int \n^\t_{a-\varepsilon}(d\sigma d\t')\,\varphi(\sigma)\,\un_{\{h(\t')\geq
\varepsilon\}}
\end{eqnarray}
\end{description}
Furthermore, for every $a>0$, the conditional distribution of the point measure
$\n^\t_a(d\sigma d\t')$, under the probability measure $\Theta(d\t\mid h(\t)>a)$ and 
given ${\rm tr}_a(\t)$,
is that of a Poisson point measure
on $\t(a)\times \T$ with intensity $\ell^a(d\sigma)\Theta(d\t')$.
\end{theorem}

The last property is the most important one. It will be called the {\it branching property}
of the L\'evy tree as it is exactly analogous to the classical branching property for
Galton-Watson trees (cf Property (ii) in the definition of Galton-Watson trees in
subsection 4.1). The random measure $\ell^a$ will be called the {\it local time}
of $\t$ at level $a$.

\smallskip
\noindent{\bf Remark.} The reader may be a little puzzled by the mathematical meaning of
the branching property as stated in the theorem, since our trees $\t$ are defined as equivalent
classes of isometric objects, and $\t(a)$ does not seem to be a well-defined object. In the 
proof below, we will circumvent this
difficulty by dealing with the tree $\t_H$ under $N$. A more
intrinsinc way to state the branching property in a mathematically precise way is as follows.
Consider first a fixed real tree $\t$, and assume that the local time measure $\ell^a$
of $\t$ can be defined via formula (\ref{appLTlevel}). Then let 
$$\sum_{j\in \jj}\delta_{(\eta_j,\t^{(j)})}$$
be a Poisson point measure on $\t(a)\times \T$ with intensity $\ell^a(d\eta)\Theta(d\t')$. 
Construct another real tree $\ov\t$ as the disjoint union
$$\ov\t={\rm tr}_a(\t)\bigsqcup\Big(\bigsqcup_{j\in\jj} (\t^{(j)}\backslash \{\rho(\t^{(j)})\})\Big)$$
equipped with the obvious appropriate distance so that 
the sets $(\t^{(j)}\backslash \{\rho(\t^{(j)}\})\sqcup\{\eta_j\}$ become the subtrees of
$\ov\t$ originating from level $a$. Note that the distribution of $\ov\t$ only depends on
the equivalence class of $\t$ in $\T$, and that this distribution is a measurable function of
$\t$. The branching property can be restated by saying that if $\t$ is
chosen according to the distribution $\Theta(d\t\mid h(\t)>a)$, then $\ov\t$ has the same
distribution as $\t$.

\medskip

\noindent{\bf Proof.} From the definition of the measure $\Theta$, it is enough to
construct the measures $\ell^a$ and to verify the properties stated in the theorem,
for the tree $\t_H$ associated with the height process $H$ under the measure $N$.  
For every $a>0$, we
define the finite measure $\ell^a$ on $\t_H$ by setting
$$\langle \ell^a,\varphi\rangle =\int_0^\zeta d\Lambda^a_s\,\varphi(p_H(s)).$$
In other words, $\ell^a$ is the image of the measure $d\Lambda^a_s$ under the mapping
$p_H$. The support property of local times guarantees that $\ell^a$ is $N$ a.e.
supported on $p_H(\{s:H_s=a\})=\t_H(a)$.

Similarly we have seen that $\{\Lambda^a_\zeta>0\}=\{\sup H>a\}$, $N$ a.e., which gives 
property (ii). Before proving property (iii), we will discuss the branching 
property of the L\'evy tree, which is basically a consequence of Corollary \ref{sublevel}.

Recall the notation $\wt H^a$ for the process $H$ truncated at level $a$, and note
the easy identification
\begin{equation}
\label{truncident}
\t_{\tilde H^a}={\rm tr}_a(\t_H).
\end{equation}
Indeed, if $A^a_s:=\int_0^s dr\,\un_{\{H_r\leq a\}}$, there is a (unique) isometry mapping
${\rm tr}_a(\t_H)$ onto $\t_{\tilde H^a}$ such that, for every $s\geq 0$ with $H_s\leq a$,
$p_H(s)$ is mapped to $p_{\tilde H^a}(A^a_s)$.

Recall the notation $\wt \h^a$ for the $\sigma$-field generated by $\wt H^a$ augmented with the
class of $N$-negligible sets.

\medskip
\noindent{\bf Claim.} {\it The conditional distribution of the point measure $\n_a^{\t_H}$, 
under the probability measure
$N(\cdot\mid\sup H>a)$ and given the $\sigma$-field $\wt \h^a$, is that of a Poisson point measure 
on $\t_H(a)\times \T$ with intensity $\ell^a(d\sigma)\Theta(d\t)$}.

\medskip

The claim is very close to the branching property (for the tree $\t_H$)
as stated in the theorem, with the minor difference that we are conditioning
with respect to the $\sigma$-field $\wt \h^a$ rather than with respect
to the truncated tree $\t_{\tilde H^a}$ (which contains less
information). We will see later that $\ell^a$ is $N$ a.e. equal to a measurable
function of the tree $\t_{\tilde H^a}$, and then the branching property
will follow from the claim.

Let us first prove the claim. We recall the notation of Corollary \ref{sublevel}: $(\alpha_j,\beta_j)$,
$j\in \jj$ are the connected  components of the open set $\{s:H_s>a\}$ and, for every $j\in \jj$,
$H^{j}_s=H_{(\alpha_j+s)\wedge \beta_j}-a$. We also set $l_j=\Lambda^a_{\alpha_j}$.
Corollary \ref{sublevel} asserts that, conditionally on $\wt\h^a$, the point measure
$$\sum_{j\in \jj} \delta_{(l_j,H^{j})}(d\ell\,d\omega)$$
is Poisson 
on $\R_+\times C_+([0,\infty))$ with intensity $\un_{[0,\Lambda^a_\zeta)}(\ell)\,d\ell\,\Delta(d\omega)$. On
the other hand, it is immediate from the construction of the tree $\t_H$ that
$$\n_a^{\t_H}=\sum_{j\in \jj} \delta_{(p_H(\alpha_j),\t_{H^{j}})}.$$
To complete the proof of the claim, we will argue that conditionally given $\wt H^a$, 
each pair $(p_H(\alpha_j),\t_{H^{j}})$ is a function (not depending on $j$) of 
the pair $(l_j,H^{j})$. This is obvious for the second coordinate $\t_{H^j}$. Then,
recalling the identification (\ref{truncident}), we have for every $j\in \jj$,
$$p_H(\alpha_j)=p_{\tilde H^a}(A^a_{\alpha_j}).$$
From (\ref{techsub}), we get for every $j\in\jj$,
\begin{equation}
\label{techbranpro}
p_H(\alpha_j)=p_{\tilde H^a}(\inf\{s:\wt \Lambda^a_s>l_j\}).
\end{equation}
This is the formula we were aiming at. In view of applying Corollary \ref{sublevel}, we still need to
determine the image
of the measure $\un_{[0,\Lambda^a_\zeta)}(l)dl$ under the mapping $l\to p_{\tilde H^a}(\inf\{s:\wt
\Lambda^a_s>l\})$.
Write $\gamma^a_s=\inf\{r:\Lambda^a_r>s\}$ and 
$\wt\gamma^a_s=\inf\{r:\wt\Lambda^a_r>s\}$, for every $s\in[0,\Lambda^a_\zeta)$. From the relation
between $\Lambda^a_s$ and $\wt\Lambda^a_s$, we see that $\wt\gamma^a_s=A^a_{\gamma^a_s}$ for
every $s\in[0,\Lambda^a_\zeta)$. Hence, via the identification (\ref{truncident}), we have also
$p_H(\gamma^a_s)=p_{\tilde H^a}(\wt\gamma^a_s)$ for every $s\in[0,\Lambda^a_\zeta)$. Therefore, for any
nonnegative measurable function $\varphi$ on $\t_H$,
$$\langle\ell^a,\varphi\rangle=\int d\Lambda^a_s\,\varphi(p_H(s))
=\int_{[0,\Lambda^a_\zeta)}dr\,\varphi(p_H(\gamma^a_r))=\int_{[0,\Lambda^a_\zeta)} dr\,\varphi(p_{\tilde
H^a}(\wt\gamma^a_s)).$$ 
Thus $\ell^a$ is the image
of the measure $\un_{[0,\Lambda^a_\zeta)}(l)dl$ under the mapping $l\to p_{\tilde H^a}(\inf\{s:\wt
\Lambda^a_s>l\})$. Using (\ref{techbranpro}) and recalling that $\Theta$ is the image of 
$\Delta(dh)$ under the mapping $h\la \t_h$, we see that
the claim follows from Corollary \ref{sublevel}.

\smallskip
We now turn to the proof of (iii). Consider first the case $\varphi=1$, where we
have to prove
$$\langle\ell^a,1\rangle =\lim_{\varepsilon\da 0}  {1\over v(\varepsilon)}\,Z(a-\varepsilon,\varepsilon)
=\lim_{\varepsilon\da 0}  {1\over v(\varepsilon)}\,Z(a,\varepsilon),$$
where $Z$ refers to the tree $\t_H$. This easily follows from the preceding claim:
To get the first equality, note that conditionally on $\langle \ell^{a-\varepsilon},1\rangle$,
$Z(a-\varepsilon,\varepsilon)$ is Poisson with parameter $v(\varepsilon)\langle
\ell^{a-\varepsilon},1\rangle$, and recall from subsection 3.2 that $\langle \ell^{a-\varepsilon},1\rangle
=\Lambda^{a-\varepsilon}_\zeta$ converges $N$ a.e. to $\langle \ell^a,1\rangle=\Lambda^a_\zeta$.
Standard estimates for Poisson variables, together with the (obvious) monotonicity of
the mapping $\varepsilon\to Z(a-\varepsilon,\varepsilon)$ then give the desired result. The case
of the second equality is treated in a similar way.

Consider then a Lipschitz function $\varphi$ on $\t_H$, with Lipschitz constant $K$. Let
$\delta>0$ be a rational number in $(0,a)$. Write $\t^{(l)}$, $l\in\l$
for the subtrees of $\t_H$ originating from level $a-\delta$, and $(\alpha_l,\beta_l)$
for the excursion interval of $H$ above level $a-\delta$ that corresponds to $\t^{(l)}$.
Again thanks to the claim, we can apply the case $\varphi=1$ of (iii) to each 
tree $\t^{(l)}$, and we get that
for every $l\in \l$,
$${1\over v(\varepsilon)}\,\n^{\t^{(l)}}_{\delta-\varepsilon}(\{(\sigma',\t'):h(\t')>\varepsilon\})
\build{\la}_{\varepsilon\to 0}^{} \Lambda^a_{\beta_l}-\Lambda^a_{\alpha_l}.$$
We have then
\begin{eqnarray*}
&&\liminf_{\varepsilon\to 0} {1\over v(\varepsilon)}
\int \n^{\t_H}_{a-\varepsilon}(d\sigma'd\t')\,\varphi(\rho(\t'))\,\un_{\{h(\t')\geq \varepsilon\}}\\
&&\qquad= \liminf_{\varepsilon\to 0} {1\over v(\varepsilon)} \sum_{l\in \l}
\int \n^{\t^{(l)}}_{\delta-\varepsilon}(d\sigma'd\t')\,\varphi(\rho(\t'))\,\un_{\{h(\t')\geq
\varepsilon\}}\\ 
&&\qquad\geq \liminf_{\varepsilon\to 0} \sum_{l\in \l} \Big({1\over v(\varepsilon)}
\n^{\t^{(l)}}_{\delta-\varepsilon}(\{(\sigma',\t'):h(\t')\geq
\varepsilon\})\,\inf_{\t^{(l)}([0,\delta])}\varphi
\Big)\\
&&\qquad= \sum_{l\in \l} \Big((\Lambda^a_{\beta_l}-\Lambda^a_{\alpha_l})\,\inf_{s\in
[\alpha_l,\beta_l],H_s\leq a}
\varphi(p_H(s))\Big)\\
&&\qquad\geq \sum_{l\in \l} \int_{\alpha_l}^{\beta_l} d\Lambda^a_s\,\Big(\varphi(p_H(s))-2K\delta\Big)\\
&&\qquad=\langle \ell^a,\varphi\rangle -2K\delta\langle \ell^a,1\rangle.
\end{eqnarray*}
The fourth line in the previous calculation is an equality because the sum is in fact finite.
In the last inequality, we used the Lipschitz property of $\varphi$, together with the
fact that the distance between $p_H(s)$ and $p_H(s')$ is bounded above by $2\delta$
whenever there exists $l\in \l$ such that  $s,s'\in[\alpha_l,\beta_l]$ and $H_s\vee
H_{s'}\leq a$.

Since $\delta$ was arbitrary we get
$$\liminf_{\varepsilon\to 0} {1\over v(\varepsilon)}
\int \n^{\t_H}_{a-\varepsilon}(d\sigma'd\t')\,\varphi(\rho(\t'))\,\un_{\{h(\t')\geq \varepsilon\}}
\geq \langle \ell^a,\varphi\rangle.$$
The same method gives the analogous bound for the limsup behavior, and a similar
argument applies to the other part of (iii). This completes the proof of
properties (i)--(iii). 

Finally, observing that for every $\varepsilon>0$
$$\int \n^{\t_H}_{a-\varepsilon}(d\sigma'd\t')\,\varphi(\rho(\t'))\,\un_{\{h(\t')\geq \varepsilon\}}
=\int \n^{\t_{\tilde H^a}}_{a-\varepsilon}(d\sigma'd\t')\,\varphi(\rho(\t'))\,\un_{\{h(\t')\geq
\varepsilon\}}$$ 
we deduce from (iii) that $\ell^a$ coincides $N$ a.e. with a measurable
function of the truncated tree ${\rm tr}_a(\t_H)=\t_{\tilde H^a}$. Consequently in the claim
above we may condition on ${\rm tr}_a(\t_H)$ rather than on the $\sigma$-field $\wt \h^a$.
This gives the branching property for $\t_H$ and completes the proof of the theorem. \cq

\smallskip
\noindent{\bf Remark.} Although our contruction of the measure $\ell^a$
for the tree $\t_H$ makes use of the coding via the height process, part (iii)
of the theorem shows that $\ell^a$ is a function of the tree and does not
depend on the particular coding that is used.

\begin{theorem}
\label{regularLT}
We can choose a modification of the collection $(\ell^a,a\geq 0)$ in such a way
that the mapping $a\la \ell^a$ is $\Theta(d\t)$ a.e. c\` adl\` ag for the weak
topology on finite measures on $\t$. We have then
$$\inf\{a> 0:\ell^a=0\}=\sup\{a\geq 0:\ell^a\ne 0\}=h(\t)\;,\quad \Theta\hbox{ a.e.}$$
\end{theorem}

\noindent{\bf Proof.} Again, it is enough to prove this for the tree $\t_H$ under $N$. 
For every rational $q>0$, denote by $\t^{(i)}$, $i\in \ii_q$ the subtrees of $\t_H$
originating from level $q$ (the index sets $\ii_q$ are disjoint when $q$ varies). 
Denote by $(\ell^{a,(i)},a\geq 0)$ the local times of $\t^{(i)}$, which are well defined
thanks to the branching property. Using the
approximations of local time given in Theorem \ref{existLT}(iii), it is immediately checked that,
for every $b>q$, we have
\begin{equation}
\label{decompLT}
\ell^b=\sum_{i\in \ii_q} \ell^{b-q,(i)}\;,\quad N\hbox{ a.e.}
\end{equation}
Note that (\ref{decompLT}) holds for every rationals
$b,q$ with $0<q<b$ outside a single set of zero $N$-measure. In addition, by 
Lemma \ref{contiLT}, we can assume that outside the same set
of zero $N$-measure, the mapping
$$\Q\cap (0,\infty)\ni b\la \langle \ell^{b,(i)},1\rangle=\ell^{b,(i)}(\t^{(i)})$$
has a c\` adl\` ag extension to $[0,\infty)$, for every rational $q>0$
and every $i\in \ii_q$.

We then show that the limit
$$\ov\ell^a:=\lim_{\Q\ni b\da a} \ell^b$$
exists for every $a\in [0,\infty)$. When $a=0$, the result is immediate, with $\ov\ell^0=0$,
since we already know that $\langle\ell^b,1\rangle\la 0$ as $b\da 0$, $N$ a.e. 
(Lemma \ref{contiLT}). Then, if $a>0$
and if $\varphi:\t_H\la \R$ is Lipschitz with constant $K$, we get from (\ref{decompLT})
that for every
$b,b'\in(a,\infty)\cap\Q$ and $q\in(0,a)\cap\Q$
\begin{eqnarray}
\label{minorLT}
\langle \ell^b,\varphi\rangle &\geq&\sum_{i\in \ii_q} \langle \ell^{b-q,(i)},1\rangle
\;\inf_{\t^{(i)}(b-q)}\varphi\;,\\
\label{majorLT}
\langle \ell^{b'},\varphi\rangle &\leq&\sum_{i\in \ii_q} \langle \ell^{b'-q,(i)},1\rangle
\;\sup_{\t^{(i)}(b'-q)}\varphi\;.
\end{eqnarray}
Moreover, we have for every $i\in \ii_q$,
\begin{equation}
\label{LipLT}
\sup_{\t^{(i)}(b'-q)}\varphi\leq \inf_{\t^{(i)}(b-q)}\varphi +K(b+b'-2q)\;,
\end{equation}
because any point of $\t^{(i)}(u)$ is at distance $u$ from the root of $\t^{(i)}$.

From the remarks of the beginning of the proof, we know that for every $i\in \ii_q$,
$$\langle \ell^{b'-q,(i)},1\rangle-\langle\ell^{b-q,(i)},1\rangle$$
tends to $0$ as $b,b'\da a$ with $b,b'\in \Q\cap [a,\infty)$. Also note that in the sums
appearing in (\ref{minorLT}) and (\ref{majorLT}) only finitely many terms
can give a nonzero contribution, namely those for which $h(\t^{(i)})>a-q$.
Combining these facts with (\ref{minorLT}), (\ref{majorLT}) and (\ref{LipLT}) leads to
$$\limsup_{\Q\ni b,b'\da a} \Big(\langle \ell^{b'},\varphi\rangle -
\langle \ell^{b},\varphi\rangle\Big)\leq 2K(a-q)\,\sup_{x\in\Q\cap (0,\infty)}\langle \ell^x,1\rangle.
$$
Since $q$ can be taken arbitrarily close to $a$, this is enough to get the convergence
of $\langle \ell^b,\varphi\rangle$ as $b\da a$, $b\in \Q$, thus proving the existence of
the right limit $\ov\ell^a$.

A similar argument gives the existence of left limits along rationals. The process 
$(\ov\ell^a,a\in [0,\infty))$ is thus c\` adl\` ag. In addition, it is easy to
see that $\ov\ell^a=\ell^a$ $N$ a.e., for every $a\in[0,\infty)$: If $a>0$ is fixed, note
that (\ref{majorLT}) will hold with $b'$ replaced by $a$ (outside a set of
zero $N$-measure depending on $a$) and use the preceding argument to verify that
the measures $\ell^a$ and $\ov\ell^a$ coincide $N$ a.e.

It remains to prove that
\begin{equation}
\label{extinctLT}
\inf\{a> 0:\ov\ell^a=0\}=\sup\{a\geq 0:\ov\ell^a\ne 0\}=h(\t_H)\;,
\quad N\hbox{ a.e.}
\end{equation}
We know that for every $\delta>0$ the (c\` adl\` ag) process $(\langle\ov\ell^{\delta+r},1\rangle,r\geq 0)$
is distributed under $N(\cdot\mid \sup H>\delta)$ as a $\psi$-CSBP. The strong Markov property
of the $\psi$-CSBP then shows that if $T=\inf\{a\geq \delta:\ov\ell^a=0\}$ we have 
$\ov\ell^b=0$ for every $b\geq a$, $N$ a.e. This gives the first equality in 
(\ref{extinctLT}). The second one is immediate from the fact that
$\{\sup H>q\}=\{\Lambda^q_\zeta\ne 0\}$, $N$ a.e., for every rational $q$. \cq

\medskip
\noindent{\bf Remark.} We already noticed that the $\psi$-CSBP has only positive jumps, and the
same holds for the total local time process $\langle\ell^a,1\rangle$ under $\Theta$. A careful inspection
of the previous proof then shows that $\Theta$ a.e. for every jump time $b$
of the process $a\longrightarrow \ell^a$ we must have $\ell^b\geq \ell^{b-}$. As
a consequence $b$ is a jump time of the process $a\longrightarrow \ell^a$ if and only
is it is a jump time of the total mass process $a\longrightarrow \langle\ell^a,1\rangle$.
More information about these jumps will be given in Theorem \ref{infinite-mult}.

\smallskip
From now on we consider only the c\` adl\` ag modification of the collection
$(\ell^a,a\geq 0)$ obtained in Theorem \ref{regularLT}. By combining the right-continuity
of the mapping $a\to \ell^a$ with Theorem \ref{existLT} (i), we get that 
$\Theta$ a.e. for every $a>0$, $\ell^a$ is supported on the level set $\t(a)$. A more
precise result will be derived below (Theorem \ref{localextinct}).

We put
\begin{equation}
\label{uniformmeasure}
{\bf m}=\int_0^\infty da\,\ell^a
\end{equation}
which defines a finite measure on the tree $\t$. When $\t=\t_H$ is the
tree coded by the height process under $N$, the measure ${\bf m}$
coincides with the image of Lebesgue measure on $[0,\zeta]$ under the
mapping $s\to p_H(s)$ (see formula (32) in \cite{DuLG}). However, formula (\ref{uniformmeasure}) 
makes it clear that the measure $\bf m$ only depends on the tree $\t$
and not on a particular coding. We will write $\zeta$
for the total mass of $\bf m$, in agreement with the case 
of the tree $\t_H$. The next theorem will imply in particular that the
topological support of $\bf m$ is $\t$, $\Theta$ a.e. 

For every $\sigma\in\t$ and $\varepsilon>0$, denote by $B(\sigma,\varepsilon)$ the open ball
of radius $\varepsilon$ centered at $\sigma$. We say that a
vertex
$\sigma$ of $\t$ is an {\it extinction point} if there exists $\varepsilon>0$ such that
$$d(\rho(\t),\sigma)=\sup_{\tau\in B(\sigma,\varepsilon)}d(\rho(\t),\tau).$$
Note that $p_H(s)$ is an extinction time of $\t_H$ iff $s\in[0,\zeta]$ is a local maximum
of $H$. As a consequence there are at most countably many extinction points. We 
denote by $\e$ the set of all extinction levels, that is levels $a$ such that
$a=d(\rho(\t),\sigma)$ for some extinction point $\sigma$.
\begin{theorem}
\label{localextinct}
The following holds $\Theta(d\t)$ a.e.:
\begin{description}
\item{\rm(i)} For every $a\in(0,h(\t)]\backslash \e$, the topological support
of $\ell^a$ is equal to the level set $\t(a)$.
\item{\rm(ii)} For every $a\in\e$, the topological support of $\ell^a$ is $\t(a)\backslash\{\sigma_a\}$,
where $\sigma_a$ is the (unique) extinction point at level $a$.
\end{description}
\end{theorem}

\noindent{\bf Proof.} First observe that if $\sigma$ is an extinction point at level $b$, the
right-continuity of the mapping $a\to\ell^a$, together with the fact that $\supp\,\ell^a\subset\t(a)$,
for every $a\geq 0$ implies that $\ell^b(B(\sigma,\varepsilon))=0$
for some $\varepsilon>0$. Therefore $\sigma\notin {\rm supp}\,\ell^b$.

The proof of the theorem will then be complete if we can prove that $\Theta$ a.e.:
\begin{description}
\item{(i)} Extinction levels of two distinct extinction points are different.
\item{(ii)} For every $a>0$ and $\sigma\in\t(a)$ which is not an extinction point,
$\ell^a(B(\sigma,\varepsilon))>0$ for every $\varepsilon>0$.
\end{description}
Recall the notation before Theorem \ref{existLT}: For every rational $q\geq 0$, 
$$\n_q^\t=\sum_{i\in \ii_q} \delta_{(\sigma_i,\t^{(i)})}$$
is the point measure associated with subtrees originating from level $q$ (the
index sets $\ii_q$ are disjoint). Set $h_i=h(\t^{(i)})$ for every 
$i\in \ii_q$ and every rational $q\geq 0$. Then it is easy to verify that
$$\e=\{q+h_i:q\in\Q_+,i\in \ii_q\}.$$
Suppose that $\sigma$ and $\sigma'$ are two distinct extinction points such that
$d(\rho(\t),\sigma)=d(\rho(\t),\sigma')=a$. Then by choosing a rational $q<a$ and sufficiently
close to $a$, we see that there exist two distinct indices $i,i'\in \ii_q$ such that
$h(\t^{(i)})=h(\t^{(i')})$. This is impossible by the branching property and the fact
that the law of $h(\t)$ under $\Theta$ has no atoms. This proves (i).

Then, if $(\ell^{b,(i)},b\geq 0)$ denote the local times associated with
the tree $\t^{(i)}$ (again this makes sense by the branching property), we get
from  Theorem \ref{regularLT} that $\ell^{b,(i)}(\t^{(i)})>0$ for every $b\in(0,h_{i})$.
Note that these properties hold simultaneously for all rationals $q\geq 0$ and $i\in \ii_q$
outside a single set of zero $\Theta$-measure. Furthermore, from the approximation of local time
given in Theorem \ref{existLT} (iii), we see that $\Theta$ a.e. for every $a>0$
and every rational $q\in(0,a)$ we have
\begin{equation}
\label{sumLT}
\ell^a=\sum_{i\in \ii_q}\ell^{a-q,(i)}.
\end{equation}
Note that in the preceding sum only finitely many terms can be nonzero. To derive this 
identity consider first the case when $a$ is also rational and then use the 
right-continuity of the mapping $a\la\ell^a$ (Theorem \ref{regularLT}). 

Finally, let $a>0$ and $\sigma\in\t(a)$. Assume that $\sigma$ is not an extinction point. Then
for every rational $q<a$, the subtree $\t^{(i)}$ originating from level $q$ and containing 
$\sigma$ is such that $h(\t^{(i)})>a-q$. Therefore $\ell^{a-q,(i)}(\t^{(i)})>0$
and by (\ref{sumLT}), $\ell^a(\t^{(i)}\cap \t(a))=\ell^a(\t^{(i)})>0$. To complete the proof
of (ii), simply note that 
$\t^{(i)}\cap \t(a)$ is contained in $B(\sigma,\varepsilon)$
if $2(a-q)<\varepsilon$. \par\cq

\subsection{Decomposing the tree along an ancestral line}

We shall need one more subtree decomposition under $\Theta$, which is 
a consequence of Corollary \ref{localkey}. Before stating it, we introduce the relevant
notation. Let $\t\in \T$ and $\sigma\in \t$. Denote by $\t^{(j),\circ}$,
$j\in\jj$ the connected components of
the open set $\t\backslash \llbracket \rho(\t),\sigma\rrbracket$, and note that for every 
$j\in \jj$, $\sigma_j:=\sigma\wedge \tau$ does not depend on the choice of $\tau\in \t^{(j),\circ}$. 
Furthermore, $\t^{(j)}:=\t^{(j),\circ}\cup\{\sigma_j\}$ is a (compact rooted) $\R$-tree
with root $\sigma_j$. The trees $\t^{(j)}$, $j\in \jj$ can be interpreted as the subtrees
of $\t$ originating from the segment $\llbracket \rho(\t),\sigma\rrbracket$. We put
$$\m_\sigma=\sum_{j\in \jj} \delta_{(d(\rho(\t),\sigma_j),\t^{(j)})},$$
thus defining a point measure on $[0,\infty)\times \T$.

\begin{theorem}
\label{Palmdec}
For every $a>0$ and every nonnegative measurable function $\Phi$ on $[0,\infty)\times \T$,
$$\Theta\Big(\int \ell^a(d\sigma)\,\exp-\langle \m_\sigma,\Phi\rangle\Big)
=\exp\Big(-\int_0^a dt\,\psi'\Big(\Theta(1-\exp-\Phi(t,\cdot))\Big)\Big).$$
\end{theorem}

\noindent{\bf Proof.} As previously, we argue on the tree $\t_H$ under $N$, and we abuse 
notation by writing $\m_\sigma$ and $\ell^a$ for the corresponding objects attached to
the tree $\t_H$. If $s\in[0,\zeta]$, we also set
\begin{eqnarray*}
\hat H^s_t&=&H_{(s+t)\wedge \zeta}\ ,\quad t\geq 0\;,\\
\check H^s_t&=&H_{(s-t)_+}\ ,\quad t\geq 0\;.
\end{eqnarray*}
Then we observe that if $\sigma=p_H(s)$, we have 
\begin{equation}
\label{Palmtech}
\m_\sigma=\p^{\hat H^s}+\p^{\check H^s}
\end{equation}
where for any continuous function $h:[0,\infty)\la[0,\infty)$ with compact support,
the point measure $\p^h$ is defined as follows.

Let $\underline h(t)=\inf_{[0,t]}h$ and let $(\alpha_i,\beta_i)$, $i\in \ii$ be the excursion intervals
of $h-\underline h$ away from $0$ (that is, the connected components of the set
$\{h-\underline h>0\}$). For every $i\in \ii$, we set
$$h^i(t)=(h-\underline h)((\alpha_i+t)\wedge \beta_i)\ ,\quad t\geq 0$$
and 
$$\p^h=\sum_{i\in \ii} \delta_{(h(\alpha_i),\t_{h^i})}.$$
The identity (\ref{Palmtech}) is then a simple 
consequence of our definitions and the construction of the tree $\t_g$ in Section 2. 

Using (\ref{Palmtech}) and Corollary \ref{localkey}, we get
\begin{eqnarray*}
&&\Theta\Big(\int \ell^a(d\sigma)\,\exp-\langle \m_\sigma,\Phi\rangle\Big)\\
&&\ =N\Big(\int_0^\zeta d\Lambda^a_s\,\exp-\langle \p^{\hat H^s}+\p^{\check H^s},\Phi\rangle\Big)\\
&&\ =e^{-\alpha a}\int \M_a(d\mu d\nu)\Big(\int Q_\mu(dh)\,e^{-\langle \p^h,\Phi\rangle}\Big)
\Big(\int Q_\nu(dh)\,e^{-\langle \p^h,\Phi\rangle}\Big).
\end{eqnarray*}
From Lemma 4.2.4 in \cite{DuLG}, we immediately get
$$\int Q_\mu(dh)\,e^{-\langle \p^h,\Phi\rangle}
=\exp-\int \mu(dt)\,N\Big(1-\exp-\Phi(t,\t_H)\Big)=\exp-\int
\mu(dt)\,\Theta\Big(1-\exp-\Phi(t,\cdot)\Big).$$
Hence,
$$\Theta\Big(\int \ell^a(d\sigma)\,\exp-\langle \m_\sigma,\Phi\rangle\Big)
=e^{-\alpha a}\int \M_a(d\mu d\nu)\,\exp\Big(-\int
(\mu+\nu)(dt)\,\Theta\Big(1-\exp-\Phi(t,\cdot)\Big)\Big).$$
Recall from subsection 3.3 that the distribution of $\mu+\nu$ under $\M_a(d\mu d\nu)$
is the law of $\un_{[0,a]}(t)dU_t$, where $U$ is a subordinator with Laplace exponent 
$\psi'-\alpha$. Therefore,
\begin{eqnarray*}
\Theta\Big(\int \ell^a(d\sigma)\,\exp-\langle \m_\sigma,\Phi\rangle\Big)
&=&e^{-\alpha a}\,E\Big[\exp-\int_0^a dU_t\,\Theta(1-\exp-\Phi(t,\cdot))\Big]\\
&=&\exp\Big(-\int_0^a dt\,\psi'\Big(\Theta(1-\exp-\Phi(t,\cdot))\Big)\Big).
\end{eqnarray*}
\par\cq

\medskip
\noindent{\bf Remark.} Combining the genealogical structure of the tree $\t$ with an
independent spatial motion leads to a construction of superprocesses which will be
explained in Section 6 below. In this setting, Theorem \ref{Palmdec} is closely related to
the representation for the historical Palm measure of superprocesses, which appears in Section 4.1 of
\cite{DP} in the stable branching case.

\medskip
We will now give a first application of Theorem \ref{Palmdec} to 
properties of the L\'evy tree.
Recall from Section 2 that for every vertex $\sigma$ of the tree $\t$, $n(\sigma)$
denotes the multiplicity of $\sigma$, defined as the number of connected components
of $\t\backslash\{\sigma\}$. We write ${\cal L}=\{\sigma\in \t\backslash \{\rho(\t)\}:n(\sigma)=1\}$
for the set of leaves of $\t$.

\begin{theorem}
\label{multiplicities}
We have:
\begin{description}
\item{\rm(i)} For every $a>0$, $\ell^a(\t\backslash\l)=0$, $\Theta$ a.e. Hence, 
${\bf m}(\t\backslash\l)=0$, $\Theta$ a.e.
\item{\rm(ii)} $n(\sigma)\in\{1,2,3,\infty\}$ for all $\sigma\in\t$, $\Theta$ a.e.
\item{\rm(iii)} The set $\{\sigma\in \t:n(\sigma)=3\}$ of binary branching points 
is empty $\Theta$ a.e. if $\beta=0$. If $\beta>0$, the set of binary branching points is
a countable
dense subset of $\t$, $\Theta$ a.e.
\item{\rm(iv)} The set $\{\sigma\in \t:n(\sigma)=\infty\}$ of infinite branching points
is nonempty with positive $\Theta$-measure iff $\pi\ne 0$. If $\langle \pi,1\rangle=\infty$,
it is $\Theta$ a.e. a countable dense subset of $\t$. If $\langle \pi,1\rangle<\infty$
it is $\Theta$ a.e. a finite (possibly empty) subset of $\t$.
\end{description}
\end{theorem}

\noindent{\bf Proof.} We can reinterpret the result of Theorem \ref{Palmdec} in the
following way. Let $U=(U_t,t\geq 0)$ be a subordinator with Laplace exponent
$$\psi'(\lambda)=\alpha+2\beta \lambda+\int \pi(dr)\,r(1-e^{-\lambda r}).$$
Note in particular that $U$ is killed at a time $\xi$ which is exponentially 
distributed with parameter $\alpha$. Write $P$ for the probability measure under
which $U$ is defined. 

Let $a>0$, and note that the formula
$$\Sigma(A)=\Theta\Big(\int \ell^a(d\sigma)\,\un_A(\m_\sigma)\Big)$$
defines a (finite) measure on the set of point measures on $[0,a]\times\T$. Then Theorem \ref{Palmdec}
says that $\Sigma$ is the law under $P(\cdot\cap\{a<\xi\})$ of a point measure 
$\m^*(dtd\t)$ which 
conditionally given $U$ is Poisson with intensity
$$\un_{[0,a]}(t)\,dU_t\,\Theta(d\t).$$

Since the measure $dU_t$ a.s. gives no mass to the singleton $\{a\}$,
$\m^*$ can have no atom of
the form $(a,\t')$ for $\t'\in \T$. It follows that
$\Theta$ a.e., $\ell^a(d\sigma)$ a.e., the point measure $\m_\sigma$ has no atom of
this form. This means that $\ell^a$-almost every $\sigma$ is 
a leaf, and (i) follows.

Recalling that $\Theta$ is an infinite measure
and using standard properties of Poisson measures, we see that 
only two possibilities may occur for instants $t\in[0,a)$:
\begin{description}
\item{$\bullet$} Either $t$ is a time of jump of $U$ and then the point measure
$\m^*$ has infinitely many atoms of the form $(t,\t')$, $\t'\in\T$.
\item{$\bullet$} Or $t$ is not a time of jump of $U$ and then the point measure
$\m^*$ has at most one atom of the form $(t,\t')$, $\t'\in\T$, and may have one
only if $\beta>0$. 
\end{description}
Using this and the relation between $\m^*$ and $\m_\sigma$, we see that
for $\ell^a$-almost every $\sigma\in \t$, the set $\rrbracket\rho(\t),\sigma\llbracket$
only contains vertices $\tau$ such that $n(\tau)\in\{2,3,\infty\}$, and the value
$n(\tau)=3$ is only possible if $\beta>0$. This property holds
simultaneously for all rationals $a$, outside a single set of zero $\Theta$-measure.

Now let $\sigma$ be any vertex in the tree $\t$. If $n(\sigma)>1$, then 
the set of descendants
of $\sigma$ is not empty and by Theorem \ref{localextinct}
we can find a rational $a>d(\rho(\t),\sigma)$ such that this 
set has positive $\ell^a$-measure. From the
preceding property we deduce that
$n(\sigma)\in\{2,3,\infty\}$. This proves property (ii).

If $\beta=0$, then we already noticed that the value $n(\tau)=3$
is not achieved by any ancestor $\tau$ of $\sigma$,  for $\ell^a$-almost every 
$\sigma\in\t$. The same argument as in the proof of (ii) now shows that
there is no vertex $\tau$ such that $n(\tau)=3$, $\Theta$ a.e. This gives the
first part of (iii).

On the contrary, if $\beta>0$, we get from the relation between $\m_\sigma$ and
$\m^*$ that the set of vertices $\tau$ with $n(\tau)=3$ is dense in
$\rrbracket\rho(\t),\sigma\llbracket$, for $\ell^a$-almost every $\sigma\in\t$.
If we apply this to all rationals $a$, we get the second part of (iii).

It remains to prove (iv). To prove this property it is 
convenient to argue on the tree $\t_H$ under $N$. If $\pi=0$, then $H$ is distributed under $N$
as a (scaled) reflected Brownian motion with drift, and the fact that there are only
binary branching points is clear since local minima of this process are
distinct (alternatively, one can also use Theorem \ref{Palmdec} in the same
way as above). Suppose then that $\pi\ne 0$. We can then note 
that for every $s$ such that $X_s-X_{s-}>0$, $p_H(s)$ is an infinite branching point 
of $\t_H$. Indeed, if $T$ is a stopping time such that $X_T-X_{T-}>0$, and $X^{(T)}$
denotes the shifted process $X^{(T)}_s=X_{T+s}-X_T$, with associated minimum
process $I^{(T)}_s$, any excursion interval of $X^{(T)}-I^{(T)}$ away from $0$
before time $\inf\{s\geq 0:X^{(T)}_s=-(X_T-X_{T-})\}$ will correspond to an excursion of
$H$ above level $H_T$, and will be associated under $p_H$
with a connected component of $\t_H\backslash\{p_H(T)\}$
(recall the discussion at the end of subsection 2.1). Since there are infinitely many such excursions,
we see that $p_H(T)$ is an infinite branching point. When $\langle\pi,1\rangle=\infty$, the set
of all jump times of $X$ is dense in $[0,\zeta]$, $N$ a.e., and it follows that 
the set of infinite branching points is dense in $\t_H$. On the other hand, when
$\langle\pi,1\rangle<\infty$, $X$ has discrete jumps and between jumps behaves like
a (scaled) Brownian motion with drift. By analysing the behavior of the process $H$
in that case, one easily obtains that the infinite branching points of $\t_H$ 
exactly correspond to jump times of $X$, so that the number of 
infinite branching points is finite $N$ a.e. \cq

\smallskip
\noindent{\bf Remark.} By using Theorem \ref{Palmdec}
as in the preceding proof, one easily gets the following
additional property. If $\int r\,\pi(dr)=\infty$, then $\Theta$ a.e. for any vertex 
$\sigma$ of $\t\backslash\{\rho(\t)\}$, the ``ancestral line'' $\llbracket\rho(\t),\sigma\rrbracket$
contains infinitely many infinite branching points. On the other hand, if 
$\int r\,\pi(dr)<\infty$, then for every $a>0$, the ancestral line of $\ell^a$-almost every vertex  
contains finitely many infinite branching points.

\medskip
We state another theorem relating discontinuities of local time to branching points
of infinite multiplicity.

\begin{theorem}
\label{infinite-mult}
Let ${\cal I}(\t)=\{\sigma\in\t:n(\sigma)=\infty\}$. Then, $\Theta$ a.e.
$\{d(\rho(\t),\sigma):\sigma\in\ii(\t)\}$ coincides with the set of discontinuity times of the 
mapping $a\longrightarrow \ell^a$. Moreover, $\Theta$ a.e., for every 
discontinuity time $b$ of the 
mapping $a\longrightarrow \ell^a$, there exists a unique $\sigma_b\in\ii(\t)\cap\t(b)$, and we have
$$\ell^b=\ell^{b-}+\lambda_{b}\,\delta_{\sigma_b}$$
where $\lambda_{b}>0$ can be obtained via the approximation
$$\lambda_{b}=\lim_{\varepsilon\to 0} {1\over v(\varepsilon)}\,n(\sigma_b,\varepsilon),$$
if $n(\sigma_b,\varepsilon)$ denotes the number of subtrees originating from 
$\sigma_b$ with height greater than $\varepsilon$.
\end{theorem}

The number $\lambda_{b}$ may be called the local time of the infinite branching point 
$\sigma_b$. As the proof will show, if the tree $\t$ is constructed as $\t=\t_H$ under the excursion measure
$N$, then
$\ii(\t)$ exactly consists of the vertices $p_H(s)$, where $s$ varies in the set of discontinuity times of
$X$.  For any such $s$, the local time of the branching point $p_H(s)$ is just the jump $\Delta X_s$ of $X$
at
$s$.

\medskip
\proof We only sketch arguments.
We assume that $\t=\t_H$ is the tree constructed from the
height process $H$ under $N$. First suppose that $\sigma\in {\cal I}(\t)$ and
let
$b=d(\rho(\t),\sigma)$. Recall the discussion at the end
of subsection 2.1. From the connection betwen the height process and the
so-called exploration process (cf Chapter 1 of \cite{DuLG}) one easily sees that there
must exist $r>0$ with $p_H(r)=\sigma$ and $\Delta X_r>0$. In a way similar to the end of
the preceding proof, we may consider the path of $X$ between $r$ and $\inf\{t\geq r:X_t=X_{r-}\}$,
and obtain that $\ell^b\geq\ell^{b-}+(\Delta X_r)\,\delta_{\sigma}$
(recall that we already know that $\ell^b\geq \ell^{b-}$).

Conversely, suppose that $b$ is a discontinuity 
time of the mapping $a\longrightarrow\ell^a$. Let $r<b$ be a rational and write
$(\t^{(i)},i\in\ii)$ for the subtrees originating from level $r$. 
For every $i\in\ii$, denote by $(\ell^{a,(i)},a\geq 0)$ the local times
of $\t^{(i)}$. From the branching property and the
fact that the
$\psi$-CSBP has no fixed
discontinuities, we get that at most one (in fact exactly one) of the 
processes $a\longrightarrow\ell^{a,(i)}$ can be discontinuous
at $b-r$. It easily follows that there
exists $\sigma_b$ with $\ell^b=\ell^{b-}+\lambda_{b}\,\delta_{\sigma_b}$,
for some $\lambda_{b}>0$. Using the
branching property at level $r_n$, for a sequence of rationals $(r_n)$ decreasing to
$b$, we get that $b$ must be a point of infinite multiplicity, and that the
approximation formula of the theorem holds for $\lambda_{b}$. We 
leave details to the reader. \cq


\medskip
As a last application of Theorem \ref{Palmdec}, we give a remarkable invariance property
of the measure $\Theta(d\t)$ under uniform re-rooting. If $\t\in\T$ and $\sigma\in\t$, we write 
$\t^{[\sigma]}$ for the ``same'' tree $\t$ with root $\sigma$. 

\begin{proposition}
\label{re-rooting}
The law of the tree $\t^{[\sigma]}$ under the measure $\Theta(d\t)\,{{\bf m}(d\sigma)\over{\bf m}(\t)}$
coincides with $\Theta(d\t)$. 
\end{proposition}

\proof It is enough to verify that, for any nonnegative measurable functional $F$ on $\T$,
\begin{equation}
\label{rooting-tech}
\Theta\Big(\int {\bf m}(d\sigma)\,F(\t^{[\sigma]})\Big)=\Theta\Big({\bf m}(\t)\,F(\t)\Big).
\end{equation}
Recall the notation ${\cal M}_\sigma$ before Theorem \ref{Palmdec}. We can easily find an explicit
``reconstruction'' functional $\Gamma$ such that $\t=\Gamma(d(\rho(\t),\sigma),{\cal M}_\sigma)$ for every
$\t\in \T$ and $\sigma\in\t$. For this functional $\Gamma$, we have also 
$$\t^{[\sigma]}=\Gamma(d(\rho(\t),\sigma),\wt{\cal M}_\sigma),$$
provided we set
$$\wt{\cal M}_\sigma=\sum_{j\in\jj} \delta_{(d(\rho(\t),\sigma)-d(\rho(\t),\sigma_j),\t^{(j)})},$$
with the notation preceding Theorem \ref{Palmdec}. Now Theorem \ref{Palmdec} implies that
for any $a>0$ and any nonnegative measurable functional $G$,
$$\Theta\Big(\int \ell^a(d\sigma)\,G({\cal M}_\sigma)\Big)
=\Theta\Big(\int \ell^a(d\sigma)\,G(\wt{\cal M}_\sigma)\Big).$$
Apply this to $G({\cal M}_\sigma)=F(\Gamma(a,{\cal M}_\sigma))$ and then integrate with respect
to $da$ to get (\ref{rooting-tech}). \cq

\section{Fractal properties of $\t$}
\subsection{Covering numbers and box counting dimensions}

Recall that for any subset $A$ of $\t$ and any $\delta >0$, we have set
$$ \n \left( A, \delta \right) =\inf \big\{ n\geq 1 \; : \; \exists \sigma_1 , \ldots , \sigma_n \in \t
\quad {\rm s.t. } \quad  A\subset \bigcup_{i=1}^n  B(\sigma_i , \delta) \Big\} .$$
The following propositions give precise rates of growth 
for $\n (\t , \delta )$ and $\n (\t(a) , \delta )$. 

\begin{proposition}
\label{covlevel} For any $a>0$, $\Theta$-a.e. on $\{ h(\t) > a\}$, 
$$ \lim_{\delta \rightarrow 0} {\n (\t(a) , \delta ) \over v(\delta)} \; = \; \ell^a(\t). $$
\end{proposition}

\begin{proposition}
\label{covt} We have $\Theta$-a.e. for all $\delta>0 $ sufficiently small 
$$ \frac{v(2\delta )}{4\delta} \, \zeta \leq \n (\t , \delta ) \leq  
\frac{12\,v(\delta /6) }{\delta} \, \zeta . $$
\end{proposition}

We immediately deduce the following corollary.

\begin{corollary}
\label{boxdim}
For any $a>0$, $\Theta$-a.e. on $\{ h(\t) > a\}$, 
$$ \diml (\t)= 1+ \diml (\t(a)) = 1 + \liminf_{\delta \rightarrow 0} \frac{\log v(\delta )}{ \log(1/ \delta)
}
$$ and 
$$ \dimu (\t)= 1+ \dimu (\t(a)) = 1 + \limsup_{\delta \rightarrow 0} \frac{\log v(\delta )}{ \log(1/ \delta)
}.
$$ 
\end{corollary}

\noindent {\bf Proof of Proposition \ref{covlevel}.} Let $a>0$ and $r\in(0,a)$. Following Section 4, we
denote by $\t^{(i)}$, $i\in \ii$ the subtrees  of $\t$ originating from level $a-r$, and we set
$$\ii_{(r)}:=\{i\in \ii:h(\t^{(i)})\geq r\}$$
so that $|\ii_{(r)}|=Z(a-r,r)$, with the notation introduced before Theorem \ref{existLT}.
We then observe that, for every $r'>r$,
$$\t(a)\subset \bigcup_{i\in \ii_{(r)}}\Big(\t^{(i)}\cap \t(a)\Big)
\subset \bigcup_{i\in \ii_{(r)}} \bar B(\rho(\t^{(i)}),r)
\subset \bigcup_{i\in \ii_{(r)}} B(\rho(\t^{(i)}),r'),$$
so that
\begin{equation}
\label{covl1}
\n(\t(a),r')\leq Z(a-r,r).
\end{equation}

On the other hand, if $\sigma$ and $\sigma'$ are two vertices in $\t(a)$ that belong
respectively to $\t^{(i)}$ and $\t^{(i')}$ for distinct indices $i$ and $i'$, we
have $d(\sigma,\sigma')\geq 2r$. Therefore, $\sigma$ and $\sigma'$ must 
belong to distinct balls of any covering of $\t(a)$ by open
balls with radius $r$. From this observation, we get
\begin{equation}
\label{covl2}
\n(\t(a),r)\geq Z(a-r,r).
\end{equation}

By combining (\ref{covl1}) and (\ref{covl2}), we see that Proposition \ref{covlevel}
follows from the case $\varphi=1$ in Theorem \ref{existLT}(iii). \cq
\medskip

\noindent {\bf Proof of Proposition \ref{covt}.} We start by proving the following lemma.

\begin{lemma}
\label{recouvre}
$\Theta$-a.e. we can find a sequence $(D_n ; n\geq 1)$ of finite subsets of $\t$ such that 
\begin{description}
\item{\rm(i)} $\quad D_n \subset D_{n+1 } \; , \; n\geq 1$.

\item{\rm(ii)} For every $ n \geq 1$ and $\sigma \in \t$, there exists $\sigma' \in D_n $ such that
$d(\sigma ,
\sigma' ) < 3 .2^{-n}$.

\item{\rm(iii)} For every $n \geq 1$ and every distinct $\sigma , \sigma' \in D_n$ we have $d(\sigma , \sigma')
\geq 2^{-n}$.

\item{\rm(iv)} ${\displaystyle \quad \lim_{n\rightarrow \infty}  2^{-n}  {|D_n| \over v(2^{-n }) }
\; = \; \zeta }$.
\end{description}
\end{lemma}

\noindent {\bf Proof of the lemma}. For any $n\geq 1 $ set
$K_n = [2^n h(\t)] -1 $. Let $n\geq 1$ and $k\in\{0,1,\ldots,K_n\}$. Denote by
$$\t^{n,k,(i)}\ ,\quad 1\leq i\leq Z(k2^{-n},2^{-n})$$
the subtrees of $\t$ originating from level $k2^{-n}$ that hit
level $(k+1)2^{-n}$. We can now use induction on $n$
to select for every $n\geq 1$, $k\in\{0,1,\ldots,K_n\}$ and $i\in\{1,\ldots,Z(k2^{-n},2^{-n})\}$,
a vertex $\sigma^n_{k,i}\in \t^{n,k,(i)}(2^{-n})\subset \t((k+1)2^{-n})$, in such a 
way that the following holds:
\begin{description}
\item{(P)} Let $k\in\{0,1,\ldots,K_n\}$ and $i\in\{1,\ldots,Z(k2^{-n},2^{-n})\}$. If $j\leq
Z((2k+1)2^{-n-1},2^{-n-1})$ is the unique index such that
$\sigma^n_{k,i}\in \t^{n+1,2k+1,(j)}$, then $\sigma^{n+1}_{2k+1,j}=\sigma^n_{k,i}$.
\end{description}
We then set
$$D_n:=\{\sigma^n_{k,i}:0\leq k\leq K_n,1\leq i\leq Z(k2^{-n},2^{-n})\}\cup\{\rho(\t)\}.$$
Property (i) is clear from (P). To prove (ii), let $\sigma\in \t$. If 
$d(\rho(\t),\sigma)\leq 2^{-n}$ the desired result is obvious
since $\rho(\t)\in D_n$.  So suppose
that $d(\rho(\t),\sigma)>2^{-n}$
and let $k\in\{0,1,\ldots,K_n\}$ be such that $(k+1)2^{-n}\leq d(\rho(\t),\sigma)<(k+2)2^{-n}$.
Clearly the ancestor of $\sigma$ at generation $k2^{-n}$ must be the root of $\t^{n,k,(i)}$
for some $i$. Then simply write
$$d(\sigma,\sigma^n_{k,i})\leq d(\rho(\t^{n,k,(i)}),\sigma)+d(\rho(\t^{n,k,(i)}),\sigma^n_{k,i})
<2. 2^{-n}+2^{-n}=3. 2^{-n}.$$
The proof of (iii) is even simpler and is left to the reader.

  It remains to prove (iv). For any $x>0$, set
$$ A_n (x) = 2^{-n} \sum_{k=1}^{[ x2^n] } 
\Big(\frac{  Z(k2^{-n}, 2^{-n}) }{ v(2^{-n})} - \ell^{k2^{-n}}(\t)\Big) . $$
Let $k\geq 1$. We know from Theorem \ref{existLT} that conditionally on the 
truncated tree ${\rm tr}_{k2^{-n}}(\t)$, 
$Z(k2^{-n},2^{-n})$ is Poisson with mean $v(2^{-n})\ell^{k2^{-n}}(\t)$. In particular, 
\begin{equation}
\label{limun}
\Theta\left( {Z(k2^{-n}, 2^{-n}) \over v(2^{-n})} -\ell^{k2^{-n}}(\t)
\;\Big|\; {\rm tr}_{k2^{-n}}(\t) \right) =0
\end{equation}
and 
\begin{equation}
\label{limdeux}
\Theta\left( \left|   {Z(k2^{-n}, 2^{-n}) \over v(2^{-n})} -\ell^{k2^{-n}}(\t) \right|^2 \;
\Big|\; {\rm tr}_{k2^{-n}}(\t)\right) = {\ell^{k2^{-n}}(\t) \over v(2^{-n})}.
\end{equation}
(In both cases the conditional expectation should be understood with respect to
the probability measure $\Theta(\cdot\mid h(\t)>k2^{-n})$.)
It is also immediate to see that $ Z(k2^{-n}, 2^{-n}) $ is a measurable
function of ${\rm tr}_{(k+1)2^{-n}}(\t)$. It easily follows
that for any $ k'>k $,
\begin{equation}
\label{limundeux}
\Theta\left(  \left( {Z(k2^{-n}, 2^{-n}) \over v(2^{-n})} -\ell^{k2^{-n}}(\t) \right) 
\left(  {Z(k'2^{-n}, 2^{-n}) \over v(2^{-n})} -\ell^{k'2^{-n}}(\t) \right) 
\right) =0.
\end{equation}
The combination of (\ref{limundeux}) and (\ref{limdeux}) gives
$$ \Theta\left( A_n (x)^2 \right) = \frac{2^{-2n}}{v(2^{-n})} 
\sum_{k=1}^{[ x2^n ]} \Theta\left( \ell^{k2^{-n}}(\t) \right) \leq \frac{x2^{-n}}{v(2^{-n})} ,$$
since $\Theta(\ell^a(\t))=N(\Lambda^a_\zeta)=e^{-\alpha a}\leq 1$, for every $a>0$ (cf the end of
Section 3).
Clearly the preceding estimate implies that $\sum_{n\geq 0} \Theta(A_n (x)^2 ) < \infty $ 
and thus, for any $x>0 $, 
\begin{equation}
\label{limps}
\Theta-{\rm a.e.} \quad \lim_{n \rightarrow \infty} A_n (x) =0.
\end{equation}
Since the mapping $b\to \ell^b(\t)$ is c\` adl\` ag, we have
$\Theta$-a.e.
$$\lim_{n \rightarrow \infty } 2^{-n } \sum_{k=1}^{[ x2^n ]} \ell^{k2^{-n}}(\t)= \int_0^x db\,
\ell^b(\t). $$ 
Together with (\ref{limps}), this implies that for any $x>0$,  $\Theta$-a.e. 
$$\lim_{n \rightarrow \infty } \frac{2^{-n}}{v(2^{-n})} \sum_{k=1}^{[ x2^n ] } Z(k2^{-n}, 2^{-n}) 
= \int_0^x db\, \ell^b(\t)   . $$
Since the height $h(\t)$ is finite, we can take $x=\infty$ in the preceding limit, which gives (iv)
since $\int_0^\infty db\,\ell^b(\t)=\zeta$. \cq

\smallskip
Proposition \ref{covt} is a simple consequence of Lemma \ref{recouvre}. Let $\delta\in(0,1/2)$
and $n\geq 1$ such that $2^{-n-1}\leq \delta<2^{-n}$. From property (iii)
in Lemma \ref{recouvre}, we get that if $\delta$ is sufficiently small,
$$\n(\t,\delta)\geq \n(\t,2^{-n-1})\geq |D_n|\geq {\zeta\over 2}\,{v(2^{-n})\over 2^{-n}}
\geq {\zeta\over 2}\,{v(2\delta)\over 2\delta}.$$
Similarly, if $3.2^{-n}\leq \delta<3.2^{-n+1}$, we get from property (ii)
in Lemma \ref{recouvre} that for $\delta$ sufficiently small,
$$\n(\t,\delta)\leq |D_n|\leq 2\zeta\,{v(2^{-n})\over 2^{-n}}\leq 2\zeta\,{v(\delta/6)\over \delta/6}.$$
This completes the proof of Proposition \ref{covt}. \cq

\subsection{Hausdorff and packing dimensions of subsets of $\t$}

We first recall the well-known inequalities 
\begin{equation}
\label{ineqdim}
\dimh (B)\leq \diml (B) \quad {\rm and }\quad \dimh (B) \leq \dimp (B) \leq \dimu (B),
\end{equation}
for any subset $B$ of $\t$ (see e.g. Chapter 3 of \cite{Fal}).
 
Let $E$ be a compact subset of $(0,\infty)$, and set $A=\sup E$. We assume that the Hausdorff
dimension and  upper box counting dimension of $E$  are equal and let $d(E)\in [0,1] $
be their common value. Recall also
the notation
$\t(E)=
\bigcup_{b\in E}
\t(b)$. The lower and upper indices $\gamma$ and $\eta$
were defined in the introduction above. The aim of this subsection is to prove the
following  theorem.
\begin{theorem}
\label{hauspacktE}
Assume that $\gamma >1 $. Then, $\Theta$ a.e. on $\{ h(\t) > A\} $, 
$$ \diml (\t(E)) = \dimh (\t(E)) = d(E) + \frac{1}{\eta -1} \quad\hbox{and} \quad  \dimu (\t(E)) = \dimp
(\t(E)) = d(E) + \frac{1}{\gamma -1} .$$
\end{theorem}

\noindent {\bf Proof}: We first get upper bounds for the box counting dimensions of $\t(E)$. Let $\delta >0
$. Analogously to the above, we use the notation $\n (E, \delta )$
for the minimal number of open intervals of length $2\delta$ that are needed to cover $E$.  We can find real
numbers
$a_i$,
$1\leq i\leq
\n(E,\delta)$ such that
$$E\subset \bigcup_{i=1}^{\n(E,\delta)} (a_i-\delta,a_i+\delta).$$
Observe that
\begin{equation}
\label{trivialcov}
\t(E\cap[0,3\delta))\subset \t({[0,3\delta)})\subset B(\rho(\t),3\delta).
\end{equation}
On the other hand, $E\cap [3\delta,\infty)$ is contained in the union of those intervals
$(a_i-\delta,a_i+\delta)$ for which $a_i>2\delta$.

Now let $b>2\delta$. Denote by $\t_{(j)}$, $1\leq j\leq Z(b-2\delta,\delta)$ the subtrees 
of $\t$ originating from level $b-2\delta$ that reach level $b-\delta$ and for every $j$, let 
$\sigma_{(j)}$ be the root of $\t_{(j)}$. Clearly any vertex in $\t({(b-\delta,b+\delta)})$
belongs to $\t_{(j)}$ for some index $j$ and thus lies within distance $3\delta$ from $\sigma_{(j)}$.
Consequently,
$$\n(\t({(b-\delta,b+\delta)}),3\delta)\leq Z(b-2\delta,\delta).$$
Now recall that $\Theta(Z(b-2\delta,\delta))=v(\delta)\Theta(\ell^{b-2\delta}(\t))\leq v(\delta)$. Hence,
$$\Theta(\n(\t({(b-\delta,b+\delta)}),3\delta))\leq v(\delta).$$
We apply this to $b=a_i$ for all indices $i$ such that $a_i>2\delta$. By summing over $i$, we get
\begin{equation}
\label{majcov}
\Theta(\n(\t({E\cap[3\delta,\infty)}),3\delta))\leq v(\delta)\,\n(E,\delta).
\end{equation}

At this point, we need the following lemma.
\begin{lemma}
\label{majv}
Assume that $\gamma >1$. Then, 
\begin{description}
\item{\rm(i)} $\displaystyle{ \quad \limsup_{\delta \rightarrow 0} {\log v(\delta ) \over\log (1/\delta )} \;
\leq 
\; {1 \over\gamma -1}}\;;
$
\item{\rm(ii)} $\displaystyle{  \quad \liminf_{\delta \rightarrow 0} {\log v(\delta ) \over\log (1/\delta ) }\;
\leq 
\; {1 \over \eta -1}}\;. $
\end{description}
\end{lemma}
\noindent {\bf Proof of the lemma.} Assertion (i) is easy from the definition of $v$ and $\gamma $. Let us
prove (ii).  If $\eta = \gamma $, (ii) is a trivial consequence of (i). So, we assume that $ \gamma < \eta $.
Let 
$ \eta' \in (\gamma , \eta ) $ and 
$\gamma' \in (1, \gamma )$. There exists a sequence $u_n \uparrow \infty $ such that $\psi (u_n) \geq
u_n^{\eta'} $ , $ n\geq 1 $.  Moreover,  for all sufficiently large $u$, we have $ \psi (u) \geq 
u^{\gamma '}$. Since $\psi $ is convex, we get for $n$ large 
enough and for any $u\geq u_n$,
$$ \psi (u) \geq  \max \left( \frac{u}{u_n} u_n^{\eta'} , u^{\gamma '} \right) .$$
Set $F(a)= \int_a^{\infty} du / \psi (u)$. The previous inequality gives, for $n$ large,
\begin{eqnarray*}
F(u_n) 
& \leq & \int_{u_n}^{\infty} \left( \max ( uu_n^{\eta '-1}, u^{\gamma '} ) \right)^{-1} du \\
 &=&  \int_{u_n}^{u_n^{\frac{\eta ' -1}{\gamma ' -1} } } \frac{du}{u u_n^{\eta '-1} } 
+\int_{u_n^{\frac{\eta ' -1}{\gamma ' -1} }}^{\infty} \;
\frac{du}{u^{\gamma '}} \\
 &\leq & C \left( u_n^{1-\eta '} \log u_n + u_n^{1-\eta '} \right)
\end{eqnarray*}
for some positive constant $C$. Hence, 
$$ \left( \liminf_{\delta \rightarrow 0} \frac{\log v(\delta )}{\log 1/\delta }  \right)^{-1} 
= \limsup_{a\rightarrow \infty } \frac{\log 1/F(a)}{\log a} 
\geq \eta ' -1 $$
and (ii) follows by letting $\eta '$ go to $\eta $. \cq

\vspace{5mm}

We deduce from the previous lemma that for any $\varepsilon >0 $,  
\begin{equation}
\label{controv}
\liminf_{\delta \rightarrow 0} \delta^{\varepsilon + \frac{1}{\eta -1 }} v(\delta ) = \limsup_{\delta \rightarrow
0} \delta^{\varepsilon + \frac{1}{\gamma -1}} v(\delta ) =0.
\end{equation}
Since $d(E)= \dimu (E)$, we also know that $\delta^{ d(E) + \varepsilon } \n (E, \delta
)$ tends to $0$ as $\delta\to 0$. Thus, if $(\delta_n)$
is any sequence of positive reals decreasing to $0$, it follows from (\ref{majcov}), (\ref{controv})
and Fatou's lemma that
$$\Theta\left( \liminf_{n\to\infty} \delta_n^{2\varepsilon + d(E) + \frac{1}{\eta -1}}  
\n (\t({E\cap[3\delta_n,\infty)}), 3\delta_n ) \right) =0.$$
Hence,
$$\liminf_{n\to\infty} \delta_n^{2\varepsilon + d(E) + \frac{1}{\eta -1}}  
\n (\t({E\cap[3\delta_n,\infty)}), 3\delta_n )=0\ ,\quad \Theta\hbox{ a.e.}$$
From (\ref{trivialcov}) we have $\n(\t(E),3\delta_n)\leq 1+\n(\t({E\cap[3\delta_n,\infty)}),3\delta_n)$
and so we get $\diml(\t(E))\leq d(E) + 1/(\eta-1) + 2\varepsilon$. Since $\varepsilon$
was arbitrary we conclude that
$\diml (\t(E)) \leq d(E) + 1/(\eta -1)$, $\Theta$ a.e.

To obtain an analogous upper bound for $\dimu (\t(E))$, we set $\delta_n = 2^{-n}$ and deduce from
(\ref{majcov}) and (\ref{controv}) that 
$$ \Theta\Big( \sum_{n\geq 1} \delta_n^{3\varepsilon + d(E) + 
\frac{1}{\gamma -1}} \n \left( \t({E\cap[3\delta_n,\infty)}) , 3\delta_n \right) \Big) < \infty. $$ 
Hence,
$$\lim_{n\to\infty}\delta_n^{3\varepsilon + d(E) + 
\frac{1}{\gamma -1}} \n \left( \t({E\cap[3\delta_n,\infty)}) , 3\delta_n \right)=0\ ,\quad\Theta\hbox{
a.e.}$$ and the bound $\n(\t(E),3\delta_n)\leq 1+\n(\t({E\cap[3\delta_n,\infty)}),3\delta_n)$ allows us to
replace $\t({E\cap[3\delta_n,\infty)})$ with $\t(E)$.
Then, a simple monotonicity argument implies that
$$ \lim_{\delta \rightarrow 0} \delta^{3\varepsilon + d(E) + \frac{1}{\gamma -1}} \n (\t(E), \delta) =0\
,\quad\Theta\hbox{ a.e.} $$ It follows that $\dimu (\t(E) ) \leq d(E) + 1/(\gamma -1)$, $\Theta$ a.e.

\smallskip
   The proof of the theorem will be complete if we verify that for any 
$\varepsilon >0 $ we have $\Theta$ a.e. on $\{h(\t)>A\}$, 
\begin{equation}
\label{lowbound}
\dimh (\t(E)) \geq d(E) + \frac{1}{\eta-1+\varepsilon}- 2\varepsilon \quad {\rm and} 
\quad \dimp (\t(E)) \geq d(E) + \frac{1}{\gamma-1+\varepsilon}- 2\varepsilon
\end{equation}
We may assume that $\varepsilon$ is small enough so that 
$1/(\gamma-1+\varepsilon)\geq 1/(\eta-1+\varepsilon)>2\varepsilon$. 
Let us prove (\ref{lowbound}). Since $\dimh (E) >
d(E) -\varepsilon
$, Frostman's lemma (see Corollary 4.12 in \cite{Fal}) gives the existence of a non-trivial finite 
measure $\nu$ supported on $E$, such that
\begin{equation}
\label{regmeas}
\forall x \in E , \; \forall \delta \in [0,1] \; : \; {\nu}( [x-\delta , x+\delta ]) 
\leq C \delta^{d(E) - \varepsilon} 
\end{equation}
where $C$ is a positive constant independent of $x$ and $\delta$. 
Define the measure $\kappa $ on $\t$ by 
$$\kappa
(d\sigma )= 
\int {\nu}(da) \ell^a(d\sigma ).$$
Then $\kappa$ is supported on $\t(E)$. Moreover $\kappa$ is finite and non-trivial
$\Theta$ a.e. on $\{ \sup H > A \}$. We will prove that $\Theta$-a.e. on $\{h(\t)>A\}$, we have 
\begin{equation}
\label{density1}
\limsup_{\delta \rightarrow 0} 
\delta^{-d(E)+ 2\varepsilon -\frac{1}{\eta -1+\varepsilon}} 
\;\kappa\left( B(\sigma , \delta ) \right) < \infty \ ,\qquad \kappa(d\sigma)\hbox{ a.e. }
\end{equation} 
and 
\begin{equation}
\label{density2}
\liminf_{\delta \rightarrow 0} 
\delta^{-d(E)+ 2\varepsilon -\frac{1}{\gamma -1+\varepsilon}} 
\;\kappa \left( B(\sigma , \delta ) \right) < \infty 
\ ,\qquad \kappa(d\sigma)\hbox{ a.e. }
\end{equation} 
Then the lower bounds (\ref{lowbound}) will follow
from classical density results for packing and Hausdorff dimensions: See e.g.
Theorems 6.9 and 6.11 in Mattila \cite{Mat} (\cite{Mat} deals with subsets
of Euclidean space, but the arguments are easily adapted to our setting).

The proof of (\ref{density1}) and (\ref{density2}) will depend on a lower bound for
the quantities
$$ \e_{\delta , \lambda , b} := \Theta \left( \int \ell^b (d\sigma ) 
e^{-\lambda \kappa (B(\sigma , \delta ))} \right)\; , \quad  \lambda , b > 0,\ \delta\in(0,1] .$$

We will apply Theorem \ref{Palmdec} in order to get this bound.
To this end, let us first fix $b>0$ and $\sigma\in\t(b)$, and use the notation introduced before
Theorem \ref{Palmdec}: $\t^{(j)}$, $j\in \jj$ are the subtrees originating from the
ancestral line $\llbracket\rho(\t),\sigma\rrbracket$, and for every $j\in \jj$,
$\sigma_j\in\llbracket\rho(\t),\sigma\rrbracket$ is the root of $\t^{(j)}$. Also 
set $d_j=d(\rho(\t),\sigma_j)$ to simplify notation.

If $\tau\in \t^{(j)}$ for some $j\in \jj$, we have $d(\sigma,\tau)=b-d_j+d(\sigma_j,\tau)$. It follows
that
$$B(\sigma,\delta)\backslash \llbracket\rho(\t),\sigma\rrbracket
=\bigcup_{j\in \jj,d_j>b-\delta} \t^{(j)}((0,\delta+d_j-b)).$$
Notice that the union in the right side is disjoint. Also observe that,
for every fixed $a>0$, $\Theta$ a.e. the measure $\ell^a$ has no atoms. Indeed, if 
$\tau\in\t(a)$ were an atom of $\ell^a$, the branching property of the L\'evy tree would imply
that $\tau$ is not a leaf, contradicting the fact that $\ell^a$ almost every 
vertex is a leaf (Theorem \ref{multiplicities}(i)). From this we get that 
$\Theta$ a.e., $\nu(da)$ a.e. $\ell^a$ has no atoms and since the set 
$\llbracket\rho(\t),\sigma\rrbracket$ has at most one point of intersection with
each level set $\t(a)$ it follows that $\kappa(\llbracket\rho(\t),\sigma\rrbracket)=0$.
Thus,
$$\kappa(B(\sigma,\delta))=\kappa(B(\sigma,\delta)\backslash \llbracket\rho(\t),\sigma\rrbracket)
=\sum_{j\in \jj,d_j>b-\delta} \kappa(\t^{(j)}((0,\delta+d_j-b))).$$
Now, if $a\leq b-\delta$ or $a\geq b+\delta$, the support
property of $\ell^a$ implies that $\ell^a(\t^{(j)}((0,\delta+d_j-b)))=0$ for every $j\in \jj$
such that $d_j>b-\delta$. On the other hand, if $b-\delta<a<b+\delta$, the approximations
of local time easily give
$$\ell^a(\t^{(j)}((0,\delta+d_j-b)))=\un_{\{0<a-d_j<\delta+d_j-b\}}\,\langle \ell^{(j),a-d_j},1\rangle,$$
where $\ell^{(j),a}$, $a>0$ obviously denote the local time measures associated with the 
tree $\t^{(j)}$. We conclude that
$$\kappa(B(\sigma,\delta))=\sum_{j\in \jj}\un_{\{d_j>b-\delta\}}
\int \nu(da)\,\un_{\{d_j<a<2d_j+\delta-b\}}\,\langle \ell^{(j),a-d_j},1\rangle.$$
In this form, we can apply the formula of Theorem \ref{Palmdec} to get 
$$ \e_{\delta , \lambda , b}=\exp\Big(-\int_0^b \psi'(\Phi_{\delta,\lambda,b}(t))\,dt\Big),$$
where
$$\Phi_{\delta , \lambda , b}(t)=\Theta\Big(1-\exp-\lambda\int
{\nu}(da)\,\un_{\{t<a<2t+\delta-b\}}\,\langle\ell^{a-t},1\rangle\Big).$$ 

Now set $\wt \psi(\lambda)={\psi(\lambda)\over \lambda}$ and note that for every $\lambda>0$,
\begin{equation}
\label{boundpsi}
\psi'(\lambda)\leq {\psi(2\lambda)-\psi(\lambda)\over \lambda}\leq 2\wt \psi(2\lambda).
\end{equation}
If $t\in[0,b-\delta]$ we have $\Phi_{\delta , \lambda , b}(t)=0$. On the
other hand, if $t\in(b-\delta,b]$ then $(t,\delta+2t-b)\subset [b-\delta,b+\delta]$ and
$$\Phi_{\delta , \lambda , b}(t)\leq \lambda
\int_{[(b-\delta)\vee t,b+\delta]}{\nu}(da)\,\Theta(\ell^{a-t}(\t))
\leq \lambda\,{\nu}([b-\delta,b+\delta]).$$
Using this bound together with (\ref{boundpsi}) we have
$$\int_0^b \psi'(\Phi_{\delta , \lambda , b}(t))\,dt\leq 2\delta\,\wt\psi(2\lambda
{\nu}([b-\delta,b+\delta])).$$
and it follows that
\begin{equation}
\label{minmin}
\e_{\delta , \lambda ,b} \geq   1-2\delta \,\widetilde{\psi } (  2\lambda {\nu}
([b-\delta , b+ \delta ])) .
\end{equation}
This is the lower bound we were aiming at.

If $r>0$, (\ref{minmin}) gives, for every $\delta\in (0,1]$, 
\begin{eqnarray*}
\Theta\left( \int \kappa (d\sigma ) \un_{ \{ \kappa (B(\sigma , \delta )) >r \}} \right)
& \leq & \frac{e}{e-1} \Theta\left( \int \kappa (d\sigma ) 
(1-e^{-\frac{1}{r} \kappa (B(\sigma , \delta ))} )
\right) \\ 
&=& \frac{e}{e-1} \int {\nu} (db ) (e^{-\alpha b} -\e_{\delta , 1/r , b})  \\
&\leq&\frac{2e}{e-1}\, \delta \int {\nu} (db ) 
\widetilde{\psi } \left(  2{\nu} ([b-\delta , b+ \delta ]) /r \right)
\end{eqnarray*}
and by (\ref{regmeas}),
\begin{equation}
\label{minoration}
\Theta \left( \int \kappa (d\sigma ) \un_{ \{ \kappa (B(\sigma , \delta )) >r \}} \right) \leq 
 C' \delta \widetilde{\psi }(  2C \delta^{d(E) -\varepsilon } /r ) 
\end{equation}
where $C'$ is a positive constant depending on ${\nu}$.  By the definition of $\eta $, for all sufficiently
large
$\lambda >0$, $\widetilde{\psi} (\lambda) \leq \lambda^{\eta -1 + \varepsilon}$.  Then, take $r=r(\delta )= 2C
\delta^{d(E) -2\varepsilon +(1/(\eta -1+ \varepsilon))}$
 in (\ref{minoration}) to get, for all sufficiently small $\delta>0$
$$ \Theta \left( \int \kappa (d\sigma ) \un_{ \{ \kappa (B(\sigma , \delta )) >r(\delta ) \}} \right) 
\leq C''
\delta^{\varepsilon (\eta -1+ \varepsilon)} .$$ Set $\delta_n =2^{-n}$. Since $\eta >1$ we deduce from the
previous inequality that
$$  \Theta\Big( \int \kappa (d\sigma ) \sum_{ n\geq 1 }
\un_{ \{ \kappa (B(\sigma , \delta_n )) >r(\delta_n ) \}} \Big) < \infty  $$
and this yields the estimate (\ref{density1}) for the upper density of $\kappa $.

It remains to prove 
(\ref{density2}). 
By the definition of $\gamma $, there exists an increasing sequence $u_n \uparrow \infty $ such that
$\widetilde{\psi } (u_n) \leq u_n^{\gamma -1 +\varepsilon } $.  Define $\delta_n$ by $u_n =
\delta_n^{\varepsilon -(1/(\gamma -1+\varepsilon))}$ and take $r(\delta_n )= 2C \delta_n^{d(E) -2\varepsilon
+1/(\gamma -1 +\varepsilon)} $ in  (\ref{minoration}) to get 
$$ \Theta\left( \int \kappa (d\sigma ) \un_{ \{ \kappa (B(\sigma , \delta_n )) 
>r(\delta_n ) \}} \right) \leq C'''\delta_n^{\varepsilon (\gamma -1+ \varepsilon)}. $$
Applying Fatou's lemma, we get that $\Theta$-a.e. for $\kappa $-a.a. $\sigma $, 
$$ \liminf_{n\rightarrow \infty } \un_{ \{ \kappa (B(\sigma , \delta_n )) > 
r(\delta_n) \} } = 0  $$
which implies the estimate (\ref{density2}) and completes the proof of the theorem. \cq

\subsection{Further results and open problems}

In this section, we briefly discuss some extensions of the preceding results.
For simplicity we restrict our attention to Hausdorff dimensions and measures.
We start by weakening the condition $\gamma>1$ in Theorem \ref{hauspacktE}.
As in Theorem 
\ref{hauspacktE} we let $E$ be a (nonempty) compact subset of $(0,\infty)$ such that
$\dimh(E)=\dimu(E)=d(E)$, and we put $A=\sup E$.
We use the standard convention ${1\over 0}=\infty$. 

\begin{proposition}
\label{weakcond}
Suppose that for every integer $k\geq 1$,
\begin{equation}
\label{weakc}
\int_a^\infty {du\over \psi(u)}=o((\log a)^{-k})\quad\hbox{as }a\to\infty.
\end{equation}
Then, $\Theta$ a.e. on $\{h(\t)>A\}$,
$$\dimh(\t(E))=d(E) + \frac{1}{\eta -1}.$$
\end{proposition}

Indeed, the proof of the estimate (\ref{density1}) does not depend on the assumption
$\gamma>1$, and this immediately gives the lower bound $\dimh(\t(E))\geq d(E) + \frac{1}{\eta -1}$.
When $\eta=1$, there is nothing more to prove. 
When $\eta>1$, a slight modification of the proof of Lemma \ref{majv}  shows that part (ii) of
this lemma still holds under the condition (\ref{weakc}). The first part
of the proof of Theorem \ref{hauspacktE} then goes through without change.

Let us consider now the general case. From the preceding remarks, one easily gets
the following statement.

\begin{proposition}
\label{generalc}
We have $\Theta$ a.e. on $\{h(\t)>A\}$,
$$d(E) + \frac{1}{\eta -1}\leq \dimh(\t(E))\leq d(E) + \liminf_{\delta\to 0} {\log v(\delta)\over
\log(1/\delta)}.$$
In particular, if $\eta=\gamma=1$, we have $\dimh(\t(E))=\infty$, $\Theta$ a.e. on $\{h(\t)>A\}$.
\end{proposition}

This leaves open the following question. Suppose that $1=\gamma<\eta$ (and that (\ref{weakc})
does not hold). Can one compute $\dimh(\t(E))$, or simply $\dimh(\t)$ ?

Finally, let us discuss the stable case where more precise results are available. 
For any suitable function $g$, we write $\hbox{\goth H}^g$ for the associated Hausdorff measure.
The following theorem is proved in \cite{DuLG2}.

\begin{theorem}
\label{Hausmeasure}
{\rm(i)} Suppose that $\psi(u)=u^2$. Set
$$g_1(r)=r\log\log(1/r)\ ,\quad g_2(r)=r^2\log\log(1/r).$$
There exist positive constants $C_1$ and $C_2$
such that $\Theta$ a.e.
$$C_1\,\zeta\leq \hbox{\goth H}^{g_2}(\t)\leq C_2\,\zeta$$
and for every $a>0$, $\Theta$ a.e. on $\{h(\t)>a\}$,
$$C_1\,\langle \ell^a,1\rangle \leq \hbox{\goth H}^{g_1}(\t(a))\leq C_2\,\langle \ell^a,1\rangle.$$
{\rm (ii)} Suppose that $\psi(u)=u^\gamma$ for some $\gamma\in (1,2)$. For every
$s>0$ set
$$h_s(r)=r^{{\gamma\over\gamma-1}}\,(\log(1/r))^{{1\over \gamma-1}}\,
(\log\log(1/r))^s.$$ 
Then, there exists a real number $\xi$ such that, $\Theta$ a.e.,
$$\begin{array}{ll}
\hbox{\goth H}^{h_s}(\t)=\infty\qquad&\hbox{if } s>{1\over \gamma-1},\\
\noalign{\smallskip}
\hbox{\goth H}^{h_s}(\t)=0&\hbox{if }s<\xi.
\end{array}$$
\end{theorem}

The construction of superprocesses
that will be developed in the next section shows that Theorem \ref{Hausmeasure}(i) is related to the 
very precise estimates which have been obtained for the Hausdorff measure
of super-Brownian motion (see \cite{Per1}, \cite{DP}, \cite{LGP} and the
references therein).

Theorem \ref{Hausmeasure}(ii) leaves open the question of determining the correct Hausdorff
measure function for $\t$ in the stable case.

\section{Some applications to super-Brownian motion}

Denote by $M_f(\R^k)$ the set of all finite measures on $\R^k$ and by 
$C_{b+}(\R^k)$ the space of all nonnegative bounded continuous functions on $\R^k$.
We also write $(P_t)_{t\geq 0}$ for the semigroup of standard Brownian motion in $\R^k$. Note
that for every $t\geq 0$, $P_t$ maps $C_{b+}(\R^k)$ into itself.

 The super-Brownian motion
with branching mechanism $\psi$ (in short the $\psi$-super-Brownian motion) is the 
(time-homogeneous) Markov process $(Z_t,t\geq 0)$ with values in $M_f(\R^k)$ whose
transition kernels can be characterized as follows. For every $\mu\in M_f(\R^k)$ and
$\varphi\in C_{b+}(\R^k)$, 
$$E[\exp(-\langle Z_t,\varphi\rangle)\mid Z_0=\mu]
=\exp(-\langle \mu,u_t\rangle),$$
where the function $(u_t(x);t\geq 0,x\in \R^k)$ is bounded and 
continuous and is the unique nonnegative
solution of the integral equation
$$u_t(x)+\int_0^t P_{t-s}(\psi(u_s))(x)\,dx=P_t\varphi(x).$$

We will now explain how the genealogical structure given by the tree $\t$ under $\Theta$
can be combined with a spatial motion to give a construction of the $\psi$-super-Brownian
motion. To present this construction in a way suitable for applications, it is convenient
to introduce the notion of a spatial tree.

Informally, a ($k$-dimensional) spatial tree is a pair $(\t,W)$ where 
$\t\in \T$ and $W$ is a continuous mapping from 
$\t$ into $\R^k$. Since we defined $\T$ as a space of equivalence
classes of trees, we should be a little more precise at this point.
If $\t$ and $\t'$ are two (rooted compact) $\R$-trees and $W$ and $W'$
are $\R^k$-valued continuous mappings defined respectively
on $\t$ and $\t'$, we say that the pairs $(\t,W)$ and $(\t,W')$ are
equivalent if there exists a root-preserving isometry $\Phi$
from $\t$ onto $\t'$ such that $W'_{\Phi(\sigma)}=W_\sigma$
for every $\sigma\in\t$. A spatial tree is then defined as an equivalent
class for the preceding equivalence relation, and we denote by $\T_{sp}$
the space of all spatial trees. Needless to say we will often abuse notation
and identify a spatial tree with an element of the corresponding equivalent class.

We denote by $\T_{sp}$ the set of all spatial trees. Recall the notation
of subsection 2.2. We define a distance on $\T_{sp}$ by setting
$$d_{sp}((\t,W),(\t',W'))={1\over 2}\;\inf_{{\cal R}\in{\cal C}(\t,\t'),(\rho,\rho')\in{\cal R}}
\Big({\rm dis}({\cal R})+\sup_{(\sigma,\sigma')\in{\cal R}}|W_\sigma-W'_{\sigma'}|\Big),$$
where $\rho$ and $\rho'$ obviously denote the respective roots 
of $\t$ and $\t'$. It is easy to verify that $(\T_{sp},d_{sp})$ is a Polish space.

Let us fix $x\in\R^k$. Also let $\t\in\T$ be a compact rooted $\R$-tree 
with root $\varnothing$ and metric $d$. We may consider the $\R^k$-valued Gaussian process
$(Y_\sigma,\sigma\in
\t)$ whose distribution is characterized by
\begin{eqnarray*}
&&E[Y_\sigma]=x\;,\\
&&{\rm cov}(Y_\sigma,Y_{\sigma'})=d(\varnothing,\sigma\wedge \sigma')\,{\rm Id}\;,
\end{eqnarray*}
where ${\rm Id}$ denotes the $k$-dimensional identity matrix. Note that
$${\rm cov}(Y_\sigma-Y_{\sigma'},Y_\sigma-Y_{\sigma'})=d(\sigma,\sigma')\,{\rm Id}.$$
From Theorem 11.17 in \cite{LT}, we know that under the condition
\begin{equation}
\label{metricent}
\int_0^1 (\log \n(\t,\varepsilon^2))^{1/2}\;d\varepsilon<\infty,
\end{equation}
the process $(Y_\sigma,\sigma\in\t)$ has a continuous modification. We keep the notation
$Y$ for this modification. Assuming that (\ref{metricent}) holds, we denote by
$Q^x_\t$ the law on $\T_{sp}$ of $(\t,(Y_\sigma,\sigma\in\t))$.

As a consequence of Proposition \ref{covt}, condition (\ref{metricent}) holds
$\Theta(d\t)$ a.e. if we assume that
\begin{equation}
\label{ent2}
\int_0^1 (\log v(\varepsilon^2))^{1/2}\,d\varepsilon <\infty.
\end{equation}
From now on, we assume that (\ref{ent2}) holds (this is automatic if $\gamma>1$,
by Lemma \ref{majv} (i)). The definition of $Q^x_\t$ then makes sense 
$\Theta(d\t)$ a.e., and we may set
$$\N_x=\int \Theta(d\t)\,Q^x_\t,$$
which defines a $\sigma$-finite measure on $\T_{sp}$. We leave it to the reader to
verify the needed measurability properties of the mapping $\t\to Q^x_\t$.

\medskip
\noindent{\bf Remark.} As a consequence of Theorem 4.5.2 in \cite{DuLG}, a necessary
and sufficient condition for the existence of a continuous modification of
the process $(Y_\sigma,\sigma\in\t)$, for $\Theta$ a.e. $\t$, should be 
$$\int_1^\infty \Big(\int_0^t \psi(u)\,du\Big)^{-1/2} dt<\infty.$$
Note that this condition is stronger than (\ref{extinct}). The proof of Theorem 4.5.2 in \cite{DuLG}
strongly depends on connections between super-Brownian motion  and partial differential equations. 
Condition
(\ref{ent2}) will be sufficient for our purposes in the present work.

\smallskip

We can now turn to connections with superprocesses. Under the measure $\N_x$, we may
for every $a>0$ define a measure $\z_a=\z_a(\t,W)$ on $\R^k$ by setting
$$\langle \z_a,\varphi\rangle=\int_0^\zeta \ell^a(d\sigma)\,\varphi(W_\sigma).$$
The next proposition reformulates a special case of Theorem 4.2.1 in \cite{DuLG}.

\begin{proposition} 
\label{represuper}
Let $\mu\in M_f(\R^k)$ and let
$$\sum_{i\in \ii} \delta_{(\t^i,W^i)}$$
be a Poisson point measure on $\T_{sp}$ with intensity $\int \mu(dx)\,\N_x$. Then the 
process $(Z_a,a\geq 0)$ defined by
\begin{eqnarray*}
&&Z_0=\mu\;,\\
&&Z_a=\sum_{i\in \ii} \z_a(\t^i,W^i)\;,\quad a>0\;,
\end{eqnarray*}
is a $\psi$-super-Brownian motion started at $\mu$. 
\end{proposition}

In the formula for $Z_a$, only finitely many terms can be nonzero, simply because
finitely many trees in the collection $(\t^i,i\in \ii)$ are such that $h(\t^i)>a$.
From Theorem \ref{regularLT}, we see that the version of $Z$ defined in the
proposition is c\` adl\` ag on $(0,\infty)$ for the weak topology on
finite measures on $\R^k$. By the known regularity properties of superprocesses
(see e.g. the more general Theorem 2.1.3 in \cite{DP}), it
must indeed be c\` adl\` ag on $[0,\infty)$. The fact that we obtain the  ``good'' version of the
superprocess is a nice feature of our construction in contrast 
with the L\'evy snake approach of
\cite{LGLJ2} or
\cite{DuLG}, where regularity properties of the resulting measure-valued process were not
immediately apparent.

In view of Proposition \ref{represuper}, the measures $\N_x$ (or rather the distribution
under $\N_x$ of the measure-valued process $(\z_a,a\geq 0)$) are called the excursion 
measures of the $\psi$-super-Brownian motion. In the quadratic branching case, these
measures play an important role in the study of connections between superprocesses
and partial differential equations: See in particular \cite{LG1}. In the case
of a general branching branching mechanism, excursion measures are constructed via
the L\'evy snake in Chapter 4 of \cite{DuLG}, and a different approach has been 
proposed recently by Dynkin and Kuznetsov \cite{DK}.

As a simple application of the representation of Proposition \ref{represuper}, we use
Theorem \ref{localextinct} to  extend a result due to Perkins \cite{Per}
in the case of the quadratic branching mechanism.

\begin{proposition}
\label{contisupport}
Let $Z=(Z_a,a\geq 0)$ be a $\psi$-super-Brownian motion in $\R^k$, and for every $a\geq 0$
let $\s_a$ denote the topological support of the random measure $Z_a$. Then the mapping
$a\la \s_a$ is c\` adl\` ag from $(0,\infty)$ into the set of all compact subsets of
$\R^k$ equipped with the Hausdorff metric. Moreover, if $a$ is a discontinuity
time of this mapping there is a point $z_a\in \R^k$ such that $\s_{a-}=\s_{a}\cup\{z_a\}$.
\end{proposition}

\noindent{\bf Remark.} If we assume that the support of $Z_0$ is compact, it is also
easy to prove that the mapping $a\la \s_a$ is right-continuous at $a=0$ for the
Hausdorff metric.

\medskip
\noindent{\bf Proof.} We may assume that $Z$ is given by the formula of 
Proposition \ref{represuper}. For every $i\in \ii$, let $\e_i$ stand for the set of 
extinction times of $\t^i$. Then each set $\e_i$ is countable, and by arguments 
similar to the proof of Theorem \ref{localextinct} it is easy to prove that the
sets $\e_i$ are pairwise disjoint. If $a\in \e_i$, write $\sigma^i_a$ for the 
extinction point of $\t^i$ corresponding to the extinction time $a$. It now follows from
Theorem \ref{localextinct} that a.s. for every $a>0$:
\begin{description}
\item{$\bullet$} ${\displaystyle\s_a=\bigcup_{i\in \ii} \{W^i_\sigma:\sigma\in \t^i(a)\}}$\quad
if\quad ${\displaystyle a\notin\bigcup_{i\in \ii} \e_i}$\ ;
\item{$\bullet$} ${\displaystyle\s_a=
\{W^j_\sigma:\sigma\in \t^j(a)\backslash \{\sigma^j_a\}\}\cup\bigcup_{i\in \ii\backslash\{j\}}
\{W^i_\sigma:\sigma\in
\t^i(a)\}}$\quad if\quad $a\in \e_j$ for some $j\in \ii$.
\end{description}
As a straightforward consequence of these formulas, one can now verify that
the mapping $a\la \s_a$ is c\` adl\` ag on $(0,\infty)$, with
$$\s_{a-}=\bigcup_{i\in \ii} \{W^i_\sigma:\sigma\in \t^i(a)\}\quad\hbox{for every }a>0.$$
Furthermore the set of discontinuity times is contained in the union of the sets
$\e_i$ over $i\in \ii$, and if $a\in \e_j$ we have
$$\s_{a-}=\s_a \cup\{W^j_{\sigma^j_a}\}.$$
We leave details to the reader. \cq

\medskip
We can also apply Theorem \ref{infinite-mult} in connection with our
construction of superprocesses. We recover the fact that for every
discontinuity time $s$ of $Z$ there is a positive real
number $\lambda_s$ and a point $\omega_s\in\R^k$ such that
$Z_s=Z_{s-}+\lambda_s\,\delta_{\omega_s}$. Precisely, there is an
index $i\in\ii$ and an infinite branching point $\sigma_s$ of $\t^i$ 
at height $s$, such that $\omega_s=W^i(\sigma_s)$ and $\lambda_s$
is the local time of the infinite branching point $\sigma_s$. We omit
details since the preceding fact is known to hold in great generality:
See Th\'eor\` eme 7 in \cite{EKR}.

\medskip

We now proceed to investigate the Hausdorff dimension of the support of
$Z_a$. From Proposition \ref{represuper}, it is enough to consider the 
random measures $\z_a$ under $\N_x$. For every $a\geq 0$, we set 
$$R_a={\rm supp}(\z_a)$$
and, if $E$ is a subset of $\R_+$,
$$R_E=\overline{\bigcup_{a\in E}{\rm supp}(\z_a)}.$$

\begin{theorem}
\label{Haus-super}
Assume that $\gamma>1$. Let $E$ be a compact subset of $(0,\infty)$ whose Hausdorff 
dimension and upper box dimension are equal to $d(E)\in[0,1]$, and set $A=\sup E$. Then,
we have
$$ \dimh R_{E} =\diml R_{E}= \left( 2d(E) + \frac{2}{\eta -1} \right) \wedge k\;, $$
$\N_0$ a.e. on $\{\z_A\ne 0\}$.
\end{theorem}

We first state a simple continuity lemma.

\begin{lemma}
\label{Holder-spat} Let $T\in\t$ be such that
$$\limsup_{\varepsilon\to 0} {\log \n(\t,\varepsilon)\over \log 1/\varepsilon}<\infty.$$
Then $Q^x_\t$ a.s., the mapping $\sigma\to W_\sigma$ is H\"older continuous with exponent
${1\over 2}-\delta$ for any $\delta\in(0,{1\over 2})$.
\end{lemma}

\noindent{\bf Proof.} Standard chaining arguments show that, for every integer $m\geq 1$
and every $u>0$, 
\begin{eqnarray*}
Q^x_\t\Big(\sup_{d(\sigma,\sigma')<u}|W_\sigma-W_{\sigma'}|\Big)
&\leq &k\Big(u^{1/2} \sqrt{\log(1+\n(\t,2^{-2m})^2)}\\
&&+16\sum_{p=m+1}^\infty
2^{-p}\sqrt{\log(1+\n(\t,2^{-2p}))}\Big).
\end{eqnarray*}
See e.g. formula (11.6) in \cite{LT} and note that a
correct choice of the distance on $\t$ in order to apply this formula is
$d'(\sigma,\sigma')=2\sqrt{d(\sigma,\sigma')}$ (see the comments on page 320  of \cite{LT}).

From the assumption of the lemma we now get the existence of a constant $C(\t)$
such that, for every $m\geq 1$ and $u\in(0,1)$,
$$Q^x_\t\Big(\sup_{d(\sigma,\sigma')<u}|W_\sigma-W_{\sigma'}|\Big)
\leq k\,C(\t)\,(u^{1/2}m^{1/2}+m2^{-m}).$$
Choosing $m$ so that $u^{-2r}<m\leq 2u^{-2r}$ we get, for every $r\in(0,1/2)$,
$$Q^x_\t\Big(\sup_{d(\sigma,\sigma')<u}|W_\sigma-W_{\sigma'}|\Big)
\leq C'(r,k,\t)\,u^{{1\over 2}-r}.$$
An application of the Borel-Cantelli lemma now completes the proof. \cq

\medskip
\noindent{\bf Proof of Theorem \ref{Haus-super}}. From the support properties
of the measures $\ell^a$, we have $\N_x$ a.e. for every $a>0$,
$${\rm supp}\; \z_a\subset\{W_\sigma:\sigma\in \t(a)\}.$$
Therefore,
$$R_E\subset\{W_\sigma:\sigma\in\t(E)\}$$
(note that $\{W_\sigma:\sigma\in\t(E)\}$ is closed as the image of the compact 
set $\t(E)$ under the continuous mapping $\sigma\to W_\sigma$).
The upper bound
$$\diml R_{E}\leq  2d(E) + \frac{2}{\eta -1} $$
is then an immediate consequence of  Theorem \ref{hauspacktE}
and Lemma \ref{Holder-spat}. Note that the assumption of
Lemma \ref{Holder-spat} holds $\Theta(d\t)$ a.e. by
Proposition \ref{covt} and Lemma \ref{majv}.

Since the bound $\dimh R_{E} \leq\diml R_{E}$ is always true, the
proof of Theorem \ref{Haus-super} will be complete if 
we can verify that
\begin{equation}
\label{lower-H}
\dimh R_{E} \geq \left( 2d(E) + \frac{2}{\eta -1} \right)\wedge k,
\end{equation}
$\N_x$ a.e. on $\{\z_A\ne 0\}$. 
To this end, let $\varepsilon>0$ and
$$b(\varepsilon)=d(E)-2\varepsilon+{1\over \eta-1+\varepsilon}.$$
As in the proof of Theorem
\ref{hauspacktE} we can consider a finite measure ${\nu}$
supported on $E$ such that, if $\kappa$ denotes the measure 
$$\kappa(d\sigma)=\int {\nu}(da)\,\ell^a(d\sigma),$$
we have 
$$\limsup_{\delta \rightarrow 0} \kappa\left( B(\sigma , \delta ) \right) 
\delta^{-b(\varepsilon)} < \infty \ ,
\qquad \kappa(d\sigma)\hbox{ a.e. }$$
Furthermore $\kappa$ is nonzero $\Theta$ a.e. on $\{h(\t)>A\}=\{\z_A\ne 0\}$.
Notice that the measure $M$ defined as the image of $\kappa$ under 
the mapping $\sigma\to W_\sigma$ is supported on $R_E$, simply because
$$M=\int {\nu}(da)\,\z_a$$
and $\nu$ is supported on $E$. 

Now, for any positive integer $q$, set
$$F_q = \{  \sigma \in \t(E) \; :\; \forall \delta \in (0 , 1/q] \; , 
\; \kappa (B(\sigma , \delta )) \leq q
\delta^{b(\varepsilon )} \} $$ 
On the event $\{h(\t)>A\}$, we can find $q_0 $ such that $\kappa (F_{q_0})>0$. 
We denote by $\widetilde{\kappa}$
the restriction of $\kappa$ to $F_{q_0}$. It is then immediate to verify that,
for any $b<b(\varepsilon)$,
\begin{equation}
\label{frost}
\int \widetilde{\kappa } (d\sigma ) \int \widetilde{\kappa } (d\sigma') d(\sigma , \sigma ')^{-b} < \infty .
\end{equation}
Finally, if $\wt M$ denotes the image of $\wt \kappa$ under
$\sigma\to W_\sigma$,  we have for any $b < k/2 \wedge b(\varepsilon ) $
\begin{eqnarray*}
Q^x_\t \left(  \int \wt M(dz ) \int \wt M(dy) |z-y|^{-2b} \right) 
& = &Q^x_\t  \left(  \int \widetilde{\kappa }
(d\sigma )
\int \widetilde{\kappa } (d\sigma') |W_{\sigma } -W_{\sigma'}|^{-2b} \right) \\
& = & \int \widetilde{\kappa } (d\sigma ) \int \widetilde{\kappa } (d\sigma')
Q^x_\t \left(  |W_{\sigma } -W_{\sigma'}|^{-2b} \right) \\
& = & C \int \widetilde{\kappa } (d\sigma ) 
\int \widetilde{\kappa } (d\sigma') d(\sigma , \sigma')^{-b} ,
\end{eqnarray*}
where $C$ is a finite constant. The latter integral is finite by (\ref{frost}).
Therefore, $\Theta(d\t)$ a.e. on $\{h(\t)>A\}$ we have $Q^x_\t$ a.s.
$$ \int \wt M(dz ) \int \wt M(dy) |z-y|^{-2b} <\infty.$$
Furthermore $\wt M$ is supported on $R_E$ because $\wt M\leq M$.
Frostman's lemma now yields the desired result (\ref{lower-H}). \cq

\end{document}